\documentclass[12pt]{amsart}
\usepackage{amsfonts}
\usepackage{amsmath}
\usepackage{amssymb, esint}
\usepackage{amsthm,color}
\usepackage{enumerate}
\usepackage{mathtools}

\newcommand{\R}{{\mathbb R}}
\newcommand{\C}{{\mathbb C}}

\newcommand{\cN}{{\mathcal N}}

\def\0{{\mathbf 0}}

\newcommand{\e}{\varepsilon}
\newcommand{\vp}{\varphi}

\newcommand{\supp}{\operatorname{supp}}

\newcommand{\diam}{\operatorname{diam}}
\newcommand{\dist}{\operatorname{dist}}

\newcommand{\grad}{\operatorname{grad}}
\newcommand{\Hess}{\operatorname{Hess}}

\theoremstyle{plain}
\newtheorem{thm}{Theorem}[section]

\newtheorem{cor}[thm]{Corollary}
\newtheorem{lem}[thm]{Lemma}

\newtheorem{prop}[thm]{Proposition}

\newtheorem{defn}[thm]{Definition}
\newtheorem{ex}[thm]{Example}

\theoremstyle{rem}

\newtheorem{case}{Case}
\newtheorem{step}{Step}

\newtheorem*{claim*}{Claim}

\newtheorem{rem}[thm]{Remark}

\numberwithin{equation}{section}

\title[Constraint   maps  with  free boundaries]{Constraint   maps  with  free boundaries: \\
the obstacle  case}

\author{Alessio Figalli}
\email{alessio.figalli@math.ethz.ch}
\address{Department of Mathematics, ETH Z\"urich,  Raemistrasse 101, 8092 Z\"urich, Switzerland }

\author{Sunghan Kim}
\email{sunghan.kim@math.uu.se}
\address{Department of Mathematics, Uppsala University, S-751 06 Uppsala, Sweden}

\author{Henrik Shahgholian}
\email{henriksh@kth.se}
\address{Department of Mathematics, KTH Royal Institute of Technology, 100 44 Stockholm, Sweden}
\thanks{The paper was finalized during the research program, ``Geometric Aspects of Nonlinear Partial Differential Equations" at Institut Mittag-Leffler. H.\ Sh.\ was partially supported by Swedish Research Council. A.\ F.\ is supported by the European Research Council under the Grant Agreement No. 721675 “Regularity and Stability in Partial Differential Equations (RSPDE)”}

\begin{document}

\begin{abstract}

This paper revives a four-decade-old problem concerning regularity theory for (continuous) constraint maps with free boundaries. Dividing the map into two parts, the distance part and the projected image to the constraint, one can prove various properties for each component. As already pointed out in the literature, the distance part falls under the classical obstacle problem, which is well-studied by classical methods. A perplexing issue,  untouched in the literature, is the properties of the projected image and its higher regularity, which we show to be at most of class $C^{2,1}$. In arbitrary dimensions, we prove that the image map is globally of class $W^{3,BMO}$, and locally of class $C^{2,1}$ around the regular part of the free boundary. The issue becomes more delicate around singular points, and we resolve it in two dimensions. In the appendix, we extend some of our results to what we call leaky maps.

\end{abstract}

\maketitle
\tableofcontents




\section{Introduction}\label{sec:intro}


Our primary objective with this paper is the potential  reincarnation of some classical results for constraint maps related to the obstacle problem. 
As we will describe below, while putting the problem into a more general framework and pushing the boundaries of the existing results, we have found new interesting questions untouched by the classics. 

Our second objective is to initiate a program to  study   constraint maps for  general types  of functionals, which may give rise to other types of free boundary conditions.
 This will be pursued in our forthcoming work for the Bernoulli type problems \cite{FKS-forthcoming}.

To formulate our model problem here, we 
let $\Omega \subset \R^n $ ($n\geq 1$) be a bounded domain, and $M \subset \R^m$  ($m\geq 2$) be a domain in the target space, with a smooth boundary $\partial M$.  
We shall consider (local)  energy-minimizing maps for functionals 
 of the simplest form,\footnote{In principle, most of our results in this paper will work with some modifications for functionals with  smooth lower order terms.} namely  
\begin{equation}\label{eq:main-energy}
\int_\Omega |Du|^2\,dx, \qquad u \in W^{1,2}(\Omega; \overline M).
\end{equation}  
Constraint maps have been studied since the 1970s, 
 especially in connection with harmonic maps, see \cite{DF}. 
 It is known that for constraint maps   $u \in W^{1,2}(\Omega;\overline M)$,  there exists a closed set $S\subset\Omega$ such that $\dim_H(S) \leq n-3$ and $u\in W_{loc}^{2,p}(\Omega\setminus S)$ for any $p<\infty$; moreover, by \cite{Duz}, a constraint map solves
\begin{equation}\label{eq:main}
\Delta u =  A_u (Du,Du)\chi_{\{u\in\partial M\}}\quad\text{in }\Omega ,
\end{equation}
in the sense of distributions, where $A$ is the second fundamental form of $\partial M$. 
The singularities of constraint maps are removable under certain topological conditions on the constraints, see e.g. \cite{Fuc2} and \cite{FW}.\footnote{For instance, one may consider a constraint map from a ball to an annulus that attains the boundary values on the outer layer of the annulus. Such a constraint map exists, but it fails to be continuous simply because the topological structures are different between balls and annuli.} In the continuity set of the map, one can also study the properties of its free boundary as in \cite{Fuc}. We also refer to \cite{Fuc3} for the partial regularity theory for more general energy-functionals. 

Constraint maps can also be considered as solutions to vectorial obstacle problems. The extension of free boundary problems into the vectorial setting  regained  attention about a decade ago. In \cite{ASUW}, energy minimizing maps for $\int( |Du|^2 + 2|u|)\,dx$ are considered; the constraint here corresponds to $M =\R^m\setminus\{0\}$ and $\partial M = \{0\}$, so the boundary is zero-dimensional. An analogous vectorial extension of Bernoulli-type problems was studied in \cite{CSY} and \cite{MTV}. Around the same time, a vectorial thin obstacle problem was studied by \cite{A}, where the boundary of the constraint is of $(n-1)$-dimension. Yet the analysis for the latter was qualitatively different (and more challenging) due  to the nonlinear structure of the governing operator. 

In this paper, we are interested in the regularity of constraint maps around their free boundaries. To have free boundaries well-defined, we confine ourselves to neighborhoods where the constraint maps are continuous. It is already known from \cite{Duz,Fuc} that the constraint maps may touch the constraint only at its concave part, i.e., $\nu_u\cdot A_u(Du,Du) \leq 0$ on $\{u \in \partial M\}$, where $\nu$ is the unit normal to $\partial M$. 
In Appendix, we shall provide further generalization of our main results to continuous weak solutions of \eqref{eq:main} that may not  fully adhere to the constraint $u\in \overline M$, in the sense that the map may leak out from  $M$.

As $\partial M$ is smooth, there is a tubular neighborhood $\cN(\partial M)$ where the nearest point projection $\Pi:\cN(\partial M) \to \partial M$ is well-defined. We then decompose a constraint map $u\in W^{1,2}(\Omega;\overline M)$ as
\begin{equation}\label{eq:decom}
u \equiv (\Pi + \rho(\nu\circ \Pi))\circ u\quad\text{in }u^{-1}(\cN(\partial M)), 
\end{equation}
where $\rho$ is the distance function to $\partial M$, and $\nu$ is the inward unit normal on $\partial M$. Thanks to this decomposition, the behavior of $u$ is fully understood by studying the {\it distance part}, $\rho\circ u$, and the {\it projected image}, $\Pi\circ u$. 

The role of the distance part, $\rho\circ u$, is understood from the classical literature, see e.g., \cite{DF} and \cite{Fuc}. It is noteworthy that $\rho\circ u$ is a solution to an almost scalar obstacle problem \eqref{eq:ur-pde}, whence it is of class $C^{1,1}$. Moreover, as the free boundary of $u$ is fully characterized by that of $\rho\circ u$, the general analysis for the free boundary of $u$ (at least at the ``basic level'') follows from the classical theory for scalar obstacle problems, see \cite{Fuc}. Consequently,  one may also deduce  the structure of the singular set, see Lemma \ref{lem:FB-sing}. Nonetheless, the recent development \cite{FS, FZ} on the fine structure of the singular set or the generic regularity \cite{FR-OS} of the free boundary are widely open in this vectorial framework. 

In contrast, the regularity of the projected image, $\Pi\circ u$, has not been considered in the literature, as far as we know. In this paper we show that, regardless of how smooth the target manifold is chosen, in general the image map $\Pi\circ u$ is at most $C^{2,1}$ (equivalently, $D^3 (\Pi\circ u)$ is at most bounded); see Example \ref{ex:C21} for the optimality of this threshold.  As we shall see below, proving such an optimal regularity for $\Pi\circ u$ is a highly challenging issue, and we believe that this paper opens up various challenging problems. 

Most of our arguments (if not all)
are robust and can be generalized so that the energy-minimizing property of the constraint map is not essential, as long as the map satisfies an  Euler-Lagrange equation  similar to  \eqref{eq:main}, without any constraint; we call them ``leaky maps''. 
Moreover, one may  also  add semilinear  terms to the right hand side of equation   \eqref{eq:main}. 
We discuss  possible extensions in Appendix \ref{app:gen}.

It is noteworthy  that constraint maps can be studied in other types of free boundary problems, and the issues with the image map appear there as well. In a forthcoming paper \cite{FKS-forthcoming} we shall discuss constraint maps for Alt-Caffarelli or Alt-Phillips type functionals.

\medskip

Our paper is organized as follows. We shall discuss our main results, and provide a heuristic discussion on our approaches, in Section \ref{sec:main}. In Section \ref{sec:prelim}, we introduce the notation and terminology, and provide a concrete setting for this paper. Section \ref{sec:reg} is devoted to the study of optimal regularity of constraint maps, where Theorem \ref{thm:C11} is proved. In Section \ref{sec:higher}, we analyze the higher regularity around the regular part of the free boundary and prove Theorem \ref{thm:higher}. In Section \ref{sec:proj}, we study the optimal regularity of the projected image of  constraint maps.
 This study  is divided  into two different subsections.
In Subsection \ref{sec:proj-nd}, we conduct analysis on a basic level, which applies to any dimension $n\geq 2$, and at the end of the section we present the proofs for Theorem \ref{thm:proj-reg} -- \ref{thm:proj-sing}. In Subsection \ref{sec:proj-2d}, we specifically treat the case $n=2$, and prove Theorem \ref{thm:C21}. In Appendix \ref{app:gen}, we discuss possible generalization of our results. Appendices \ref{app:schwarz} and \ref{app:newt} are devoted to the proofs for some technical lemmas stated in Section \ref{sec:proj-2d} but proved later to maintain the reading flow in Section \ref{sec:proj-2d}. 

\section{Main results}\label{sec:main}

In all this paper, we assume $\Omega$ to be a bounded domain in $\R^n$, $n\geq 1$, and $M$ an open set in $\R^m$, $m\geq 2$, with $\partial M$ of class $C^\infty$. The smoothness of $\partial M$ can be made sharper in each statement, but this is left out for clarity of exposition.

\subsection{The constraint map}

We first prove that if $u$ is continuous in a neighborhood of a free boundary point, then $\Pi\circ u$ admits continuous second derivatives. An intuitive way to understand this is noting that the target manifold $\partial M$ flattens out as we zoom into a free boundary point. Combining this observation with the optimal regularity of $\rho\circ u$, we readily obtain the optimal regularity for any  constraint map around points of continuity.
  Although the analysis does not involve any technical difficulties, we were not able to find any literature presenting this result, so we include  it here. 

\begin{thm}\label{thm:C11} (Essentially classical)
Let $u\in W^{1,2}(\Omega;\overline M)$ be a local constraint map with $\partial M$. Then there exists a closed set $S\subset \Omega$ such that $\dim_H(S) \leq n-3$, and $u \in C^{1,1}_{loc}(\Omega\setminus S)$. 
\end{thm}

Despite the fact that $C^{1,\alpha}$-regularity of the free boundary around a point with positive density is already proved in \cite{Fuc}, the higher-regularity does not seem to be treated in the literature.
 Unlike the $C^{1,\alpha}$-regularity, which follows directly from the result for scalar obstacle problems, the higher regularity needs some bootstrap argument not needed in the scalar case. As a matter of fact, one has to combine the partial hodograph-Legendre transformation and the higher-order regularity estimates for elliptic transmission problems. Our proof however does not yield the analyticity of the free boundary, which seems to be a new challenge in the vectorial setting.

\begin{thm}\label{thm:higher}
Let $u \in W^{1,2}(\Omega;\overline M)$ be a local constraint map of \eqref{eq:main-energy}, and assume that $x_0\in (\Omega\setminus S)\cap \partial\{u\in M\}$ is regular.\footnote{See Definition \ref{def:FB}.} Then there is a ball $B\subset \Omega$ centered at $x_0$ such that $B\cap \partial\{ u\in M\}$ is a $C^\infty$-graph, and $u\in C^\infty(B\cap\overline{\{ u\in M\}})\cap C^\infty(B\cap\{ u \in\partial M\})$. 
\end{thm}

\begin{rem}\label{rem:higher}
By definition, $|Du(x_0)| > 0$ at a regular point $x_0 \in\Omega\cap \partial \{ u \in M\}$; see Definition \ref{def:FB}. Hence, by the inverse function theorem, if $m \geq n$ then there is a ball $B\subset\Omega$ centered at $x_0$ such that $u(B\cap \partial \{ u\in M\})$ is a $C^\infty$-hypersurface in $\R^m$. 
\end{rem}

\subsection{The image map in arbitrary dimensions}

We next investigate the regularity of the image map $\Pi\circ u$.
We prove that  $\Pi\circ u$ is at most $C^{2,1}$, regardless of how smooth the boundary of the target $\overline M$ is. In fact,  in Example \ref{ex:C21} we construct a constraint map that verifies this regularity threshold. As a next step, we show in any ambient space that $D^3 (\Pi\circ u)$ is of class $BMO$, and of class $L^\infty$ around a neighborhood of the regular part of the free boundary. Hence, the optimal regularity for the image map is always achieved around the regular points.\footnote{By the generic regularity techniques developed by Ros-Oton, Serra and the first author \cite{FR-OS}, it is natural to conjecture that, in dimensions less than or equal to $4$, the free boundary generically consists only of regular points and therefore the image map is class $C^{2,1}$. Such an extension of the techniques of \cite{FS,FR-OS} to this setting would be extremely interesting.
Note that we expect the restriction on the dimension to be only on the ambient space, not on dimension of the target.} We summarize this in the next statement.

\begin{thm}\label{thm:proj-reg}
Let $u \in W^{1,2}(\Omega;\overline M)$ be a local constraint map of \eqref{eq:main-energy}, and let $\Pi $ be the nearest point projection onto $\partial M$ in a tubular neighborhood $\cN(\partial  M)$. Suppose that $u\in C(B; \cN(\partial  M))$ in a ball $B\subset\Omega$. Then $D^3 (\Pi\circ u) \in BMO_{loc}(B)$. Moreover, if every point $x_0\in \partial\{ u \in M \}\cap B$ is regular, then $D^3 (\Pi\circ u)\in L_{loc}^\infty(B)$. 
\end{thm}

As a next step, we analyze the behaviour of $\Pi\circ u$ around singular points.\footnote{These are the countinuity points of $u$ where the free boundaries is not regular, see Definition \ref{def:point-reg}.} We shall observe that, if the distance part $ \rho \circ u$ is approximated by a quadratic polynomial with order $2+\bar\sigma$, then the image map $\Pi\circ u$ is approximated by a cubic polynomial with order $3+\sigma$. By the fine structure of the singular sets for the scalar obstacle problem, first investigated by Serra and the first author \cite{FS} and then refined by Franceschini and Zaton \cite{FZ}, our result suggests pointwise $C^{3,\sigma}$-regularity of $\Pi\circ u$ on the singular part of the free boundary, up to a set of Hausdorff dimension $n-3$ (this restriction comes from the presence of an exceptional set in dimensions larger than or equal to $3$, see \cite{FS}). 

\begin{thm}\label{thm:proj-sing}
Under the setting of Theorem \ref{thm:proj-reg}, suppose that $x_0\in\partial\{u\in M\}\cap B$ is a singular point.\footnote{See Definition \ref{def:FB}.} Let $\rho$ be the distance function to $\partial M$. 
If $\rho\circ u\in C^{2,\bar\sigma}(\{x_0\})$ for some $\bar\sigma \in(0,1)$, then $\Pi\circ u \in C^{3,\sigma}(\{x_0\})$ for some $\sigma\in (0,\bar\sigma)$, depending only on $n$ and $\bar\sigma$. 
\end{thm}


\subsection{The image map in two dimensions}

Our last result concerns the optimal regularity of $\Pi\circ u$ in dimension $2$. In two dimensions, neither the set $S$ of discontinuity nor the exceptional set in the singular part of the free boundary exist \cite{DF,Fuc,weiss,FS}. Thanks to Theorem~\ref{thm:proj-sing}, the problem boils down to prove a uniform $C^{2,1}$-estimate around regular points, regardless of presence of a singular point in the vicinity. 
This is a challenging  issue,
 as the local $C^{2,1}$-estimate around regular points depend on the boundary values of the image map. In fact, the behavior of regular points approaching singular points can be rough, and has been studied qualitatively only recently, by Eberle, Weiss and the third author \cite{ESW2}. 
 Since our analysis requires fine quantitative estimates, we will introduce a uniform geometric approximation property (see Definition \ref{def:gap}) for solutions to the scalar obstacle problems, which provides a quantitative estimate between the local and global solutions.
 
 Our result reads as follows.

\begin{thm}\label{thm:C21}
Let $u\in W^{1,2}(\Omega;\overline M)$ be a local constraint map of \eqref{eq:main-energy} in a bounded domain $\Omega\subset \R^2$, with a smooth open set $M\subset\R^m$. If $\rho\circ u$ admits the geometric approximation property uniformly in $\Omega$, then $D^3(\Pi\circ u) \in L_{loc}^\infty(\Omega)$.
\end{thm}

It turns out that the geometric approximation property is always true in dimension 2, thanks to a recent result of Eberle and Serra \cite{ES}. Hence, Theorem \ref{thm:C21} implies the optimal regularity of the image map in two dimensions.
\begin{cor}\label{cor:C21}
Let $u\in W^{1,2}(\Omega;\overline M)$ be a local constraint map of \eqref{eq:main-energy} in a bounded domain $\Omega\subset \R^2$, with a smooth open set $M\subset\R^m$. Then $D^3(\Pi\circ u) \in L_{loc}^\infty(\Omega)$.
\end{cor}

\subsection{Heuristic discussions on the image map}
To prove our results on the image map $\Pi\circ u$, we first write down a PDE system for it, see  \eqref{eq:ut-sys}.
Neglecting the higher-order terms in the system, the study of optimal regularity for $\Pi\circ u$ boils down to that of the solutions to the system 
\begin{equation}\label{eq:v-pde-gen}
\begin{dcases}
\Delta v  = a \cdot Dw, \\
\Delta w = g \chi_{\{w >0\}},\,\, w\geq 0, 
\end{dcases}
\quad\text{in }\Omega,
\end{equation}  
where both $a,g\in C^{1,\alpha}(\Omega)$, with $g\geq 0$.
From the regularity theory of the obstacle problem, $w \in C^{1,1}$. Hence $\Delta v \in W^{1,\infty}$,
and so elliptic regularity implies that $D^3v \in BMO.$ In other words, modulo technical details, proving Theorem~\ref{thm:proj-reg} is relatively easy, while showing the optimal $C^{2,1}$ regularity is extremely challenging. In particular, if one hopes to prove such a result, one needs to exploit finer properties of $w$.

We note that $v\in C_{loc}^{2,1}(\Omega)$ is equivalent with $v$ being approximated of order three  by quadratic polynomials, locally uniformly for all free boundary points, i.e., 
\begin{equation}\label{eq:v-qx0}
\sup_{B_r(x_0)\cap\Omega} |v - q_{x_0}| \leq \frac{Cr^3}{\dist(x_0,\partial\Omega)^3},\quad\forall r > 0, \,x_0\in\partial\{w> 0\}\cap\Omega,
\end{equation}
where $C$ does not depend on $x_0$ and $r$. 
As we shall see, although $g$ may vanish at some free boundary points of $w$, estimate \eqref{eq:v-qx0} becomes much easier at those points (Lemma \ref{lem:rescale}). In fact, it turns out to be sufficient to consider the case where $ 0< c_0 < g < c_1$ in  $\Omega$ (Remark \ref{rem:rescale}).

Now, given a free boundary point of $w$, we can distinguish whether this point is regular or singular.
If $x_0$ is singular, the fine regularity of the singular set \cite{Caff-obs,weiss,colombo--,FS,FZ} tells us that
\begin{equation}\label{eq:w-px0}
\sup_{B_r(x_0)\cap\Omega} |w - p_{x_0}| \leq Cr^2\psi(r),\quad \forall r > 0, 
\end{equation}
where $\psi$ is a universal modulus of continuity, and $p_{x_0}$ is the quadratic polynomial given by the (quadratic) blowup of $w$ at $x_0$. Also, in two dimensions, \eqref{eq:w-px0} holds with a universal H\"older modulus $\psi$. Thanks to this fact and a standard approximation argument (Lemma \ref{lem:below}),   in two dimensions we are able to prove that \eqref{eq:v-qx0} holds uniformly at every singular point of $w$.

On the other hand, if $x_0$ is regular, then $\partial\{w>0\}\cap B_{\delta_{x_0}}(x_0)$ is a $C^{1,\alpha}$-graph for some $\delta_{x_0}>0$ that may vary upon $x_0$, and so  \eqref{eq:v-qx0} holds but with a constant depending on $\delta_{x_0}$.
Now, if $\partial\{w > 0\}\cap\Omega$ does not contain any singular point (so that the entire free boundary is a single piece of $C^{1,\alpha}$-hypersurface), we can find a uniform lower bound of $\delta_{x_0}$, whence \eqref{eq:v-qx0} is verified for all $x_0\in\partial\{w>0\}$, finishing the proof.\footnote{Due to the recent breakthrough in the generic regularity \cite{FR-OS} for the (scalar) obstacle problem, the set of singular points for generic free boundaries has zero $\mathcal H^{n-4}$-measure; thus, for dimensions $n\leq 4$,  
one may expect that $v\in C_{loc}^{2,1}(\Omega)$ in the generic sense.} 
However, in general, singular points exist. So, we need to prove a bound around regular points that do not degenerate as we approach a singular point. For this, it is crucial for us to understand how the free boundary looks like at {\it every} free boundary point, and this is why we shall introduce the {\it geometric approximation property} for solutions to the scalar obstacle problems. This property can be defined for any dimension $n$.

\begin{defn}[Geometric Approximation Property]\label{def:gap}
Let $g$ be a nonnegative continuous function on a bounded domain $\Omega\subset\R^n$, and let $w$ be a nonnegative solution to $\Delta w = g \chi_{\{w > 0\}}$ in $\Omega$. The solution $w$ is said to admit the geometric approximation property at $x_0\in\partial\{w > 0\}\cap\Omega$, with order $2+\sigma$, if there is a constant $\lambda > 0$ for which 
\begin{equation}\label{eq:w-Ur-re}
\| w- U_{x_0}^r\|_{C^1(B_r(x_0)\cap\Omega)}^* \leq \lambda \max\{r^3, g(x_0)^{1-\sigma} r^{2+\sigma}\},\quad\forall r \in (0,{\rm dist}(x_0,\partial\Omega)),\footnote{See \eqref{eq:adim} for the definition of $\|\cdot\|_{C^1}^*$. }
\end{equation} 
where $U_{x_0}^r$ is a nonnegative global solution to $\Delta U_{x_0}^r = g(x_0)\chi_{\{ U_{x_0}^r > 0\}}$ in $\R^n$ such that $U_{x_0}^r (x_0) = 0$, and it may vary upon each $r>0$. 
\end{defn}

The right-hand side of \eqref{eq:w-Ur-re} is chosen to take into account the situation where $g(x_0) \ll 1$, see Remark \ref{rem:rescale}. As a matter of fact, the geometric approximation property is always true for scales $r \geq g(x_0)\geq 0$, provided that $g$ is Lipschitz in $\Omega$, see the proof of Lemma \ref{lem:rescale}. 

Our argument to prove Theorem \ref{thm:C21}, which is contained in Section~\ref{sec:proj-2d}, proceeds as follows: we fix a free boundary point $x_0$, we start from $r=1$, and we look at the size of the contact set inside $B_r(x_0)$. As long as the contact set is sufficiently small inside $B_r(x_0)$ (in other words, $x_0$ looks like a singular point at all scales $r \in (r_0,1)$), we can show that Lemma \ref{lem:below} holds up to $r=r_0$. If $r_0=0$ it means that $x_0$ is singular and we are done. If $r_0>0$, we can exploit the fact that the contact set has some thickness and the geometric approximation property to show the validity of \eqref{eq:v-qx0} for $r \in (0,r_0)$. We note that this last part, which is the most delicate, makes use of the complete classification of the global solutions in two dimensions and their explicit representation via the Schwarz functions (due to M. Sakai \cite{S}), together with some explicit estimates on the generalized Newtonian potential (Lemma \ref{lem:quad}). The proof of Theorem \ref{thm:C21} then follows after a standard approximation argument.

\begin{rem}
As mentioned before, the validity of the geometric approximation property in two dimensions has been recently obtained in \cite{ES}. Due to the presence of anomalous points in dimension $n \geq 3$ (see \cite{FS}), it is not clear whether the geometric approximation property holds for $n\geq 3$. However, one may hope it to be true if the free boundary of $w$ contains no anomalous points. We formulate this as a conjecture:

\smallskip

\noindent
{\em Conjecture:} Let $w$ be a solution to $\Delta w = g\chi_{\{w > 0\}}$ in a bounded domain $\Omega$ in $\R^n$, $n\geq 2$, for some positive Lipschitz function $g$ in $\Omega$. Denote by $\Sigma$ the singular part of $\partial\{w>0\}\cap \Omega$. Assume there exists a neighborhood $N\Subset\Omega$ and a constant $\bar\sigma\in(0,1)$ such that $w \in C^{2,\bar\sigma}(\Sigma\cap \overline N)$.\footnote{See Definition \ref{def:point-reg}.} Then $w$ verifies the geometric approximation property locally with order $2+\sigma$ uniformly in $N$, where $\sigma = \sigma (n, \bar\sigma)>0$. 

\smallskip

\noindent
It is also worth noticing that, while our proof strongly relies on the geometric approximation property, the validity of the $C^{2,1}$ regularity of the image map may be independent of it. It would be extremely interesting to understand whether this is the case or not.
\end{rem}

\subsection{Optimal regularity of the image map: an example}
We construct a constraint map of one-variable whose image map is of class $C^{2,1}\setminus C^3$ around its free boundary point. 

\begin{ex}\label{ex:C21}
Let $n = 1$, $m = 2$, $M := \R^2\setminus B_1$, so that $\partial M \equiv \partial B_1$, and consider the curve  
$$
u (x) = \begin{cases}
(1,-x) & \text{if } x < 0, \\
(\cos x, -\sin x), &\text{if } x > 0.
\end{cases}
$$
By direct computation, one can check that 
\begin{equation}\label{eq:polar}
u_{xx} =  - |u_x|^2 u\chi_{\{ |u| = 1\}}\quad\text{in }(-\infty,\infty),
\end{equation} 
and that $u$ is locally minimizing the Dirichlet energy \eqref{eq:main-energy} in any small interval around $0$. Note that $0$ is the (only) free boundary point of $u$. The image map, i.e., the nearest point projection of the constraint map to the circle, is given by 
$$
V(x) := (\cos\theta(x), - \sin\theta(x))\quad\text{with} \quad \theta(x)= 
\begin{cases} 
\tan^{-1}(x), &\text{if } x< 0, \\ 
x, &\text{if } x >0.
\end{cases}
$$
Clearly, $V_{xxx} \in L^\infty$, but the continuity of $V_{xxx} $ breaks at $x = 0$.
\end{ex}


\section{Preliminary analysis}\label{sec:prelim}

\subsection{Notation and terminologies}

Throughout this paper, $n$ denotes the dimension of the ambient space, and $m$ that of the target space, with $n\geq 1$ and $m\geq 2$. By $\min\diam E$ we denote the minimal diameter of a set $E$, i.e., the width of the narrowest strip containing $E$. Derivatives in the ambient space will be denoted by $D^k$, while $\nabla^k$ shall be used to denote those in the target space. Let $\| \cdot \|_{C^1}^*$ be the adimensional $C^1$-norm: 
\begin{equation}\label{eq:adim}
\| f \|_{C^1(B_r(x_0))}^* := \sup_{B_r(x_0)} ( |f| + r |Df|). 
\end{equation}  
Given a closed set $F$, we use $C^{k,\sigma}(F)$ to denote the class of functions with $C^{k,\sigma}$-approximation at each point on $F$: 

\begin{defn}\label{def:point-reg}
Let $\Omega\subset \R^n$ be open, let $w \in C(\Omega)$, and let $F\subset\Omega$ be a closed set. We say $w\in C^{k,\sigma}(F)$, if there exists a constant $\lambda > 0$ such that for each $x_0\in F$, one can find a polynomial $p_{x_0}$ of degree $k$ for which
$$
\sup_{B_r(x_0)\cap\Omega} |w - p_{x_0}| \leq \lambda r^{k+\sigma},\quad\forall r \in (0,1).
$$
We remark that the definition of $C^{k,\sigma}(F)$ is equivalent to the usual definition of H\"older spaces when $F$ is replaced by an open set. 
\end{defn}

Let $M \subset \R^m$ be an open set with boundary $\partial M$ of class $C^\infty$, $\rho$ the distance function to $\partial M$, $\nu$ the unit normal to $\partial M$ pointing inwards $M$, and $A$ the second fundamental form on $\partial M$. Then there exists a tubular neighborhood, $\cN(\partial  M)$, where $\rho$ is smooth. Moreover, one can define  in $\cN(\partial  M)$  the nearest point projection, $\Pi$, onto $\partial M$, which is also smooth in the tubular neighborhood; see \cite[Lemma 14.16]{GT}. 

We shall fix the tubular neighborhood $\cN(\partial  M)$ throughout the paper. With these definitions, denoting by $id$ the identity map on the target space,   
\begin{equation}\label{eq:setting}
id = \Pi + \rho(\nu\circ \Pi)\quad\text{in }\cN(\partial  M).
\end{equation}
Therefore, the decomposition in \eqref{eq:decom} follows. 

Given a map $u:\Omega\to \overline M$, to simplify the notation we set
\begin{equation}\label{eq:setting-re}
V := \Pi \circ u,\quad w := \rho\circ u\quad\text{in } u^{-1}(\cN(\partial  M)). 
\end{equation} 
We can also decompose the Jacobian matrix $Du$ as 
\begin{equation}\label{eq:Du-decom}
Du \equiv  (Du)^\tau + (Du)^\nu\quad\text{in }u^{-1}(\cN(\partial  M)),
\end{equation}
where for each $1\leq \alpha\leq n$, the vectors $(D_\alpha u)^\tau$ and $(D_\alpha u)^\nu$, evaluated at $x$, are the tangential and respectively normal component of $D_\alpha u(x)$ with respect to the tangent hyperplane of $\partial M$ at $V(x)$. More precisely, 
\begin{equation}\label{eq:Dut-Dun}
(Du)^\tau = DV + w D(\nu\circ V),\quad (Du)^\nu = (Dw) \nu_V,
\end{equation}
with $V$, $w$, and $\nu$ as in \eqref{eq:setting} and \eqref{eq:setting-re}. 

A constraint map for \eqref{eq:main-energy} is defined as follows.

\begin{defn}\label{def:map} 
Given a domain $\Omega\subset\R^n$, we shall call $u \in W^{1,2}(\Omega;\overline M)$ a local constraint map if $$\int_\Omega |Du|^2 \,dx \leq \int_\Omega |Dv|^2 \,dx,$$ for every $v\in W^{1,2}(\Omega;\overline M)$ with $\supp(u-v)\Subset \Omega$. 
\end{defn} 

\subsection{Known results}

Let us retrieve some classical results.

\begin{thm}[\cite{DF}, \cite{Fuc}]\label{thm:prelim}
Let $u\in W^{1,2}(\Omega;\overline M)$ be a local constraint map. Then there exists a closed set $S\subset\Omega$, whose Hausdorff dimension is at most $n-3$, such that $u\in W_{loc}^{2,p}(\Omega\setminus S)$ for any $p\in(1,\infty)$, and it solves \eqref{eq:main} in the sense of distributions. Moreover, for any ball $B\subset\Omega$ with $u(B)\subset \cN(\partial M)$, $w := \rho\circ u \in C_{loc}^{1,1}(B\setminus S)$ and 
\begin{equation}\label{eq:ur-pde}
\Delta w  =  - \nu_V\cdot A_V((Du)^\tau,(Du)^\tau) \chi_{\{ w > 0\}}\quad\text{in }B,
\end{equation}
in the sense of distributions. Furthermore, if $\{u \in \partial M\}$ has positive density at a point $x_0\in\partial\{u\in M\}\cap\Omega$
and $|Du(x_0)|>0$, then $\partial \{u\in M\}\cap B$ is a $C^{1,\alpha}$-graph for a ball $B\Subset\Omega$ centered at $x_0$ and some $\alpha \in(0,1)$. 
\end{thm}

\begin{rem}\label{rem:prelim}
As $w := \rho\circ u$ is nonnegative and $w = 0$ on the contact set $\{u \in\partial M\}$, due to \eqref{eq:ur-pde}, it is known (e.g., see \cite{Duz}) that 
$$
-\nu_u \cdot A_u(Du,Du) \geq 0 \quad\text{on }\{u\in \partial M\};
$$
that is, $u$ touches $\partial M$ only at the points in neighborhoods of which $\partial M$ can be represented by concave graphs; 
note that $V = u$ and $(Du)^\tau = Du$ a.e.\ on $\{u\in\partial M\}$. 

\end{rem}

In a neighbourhood where a constraint map is continuous, one can define its free boundary as follows. 

\begin{defn}\label{def:FB}
Given a ball $B$ and a map $u\in C(B;\cN(\partial  M))$, we call  $\partial\{ u\in M\}\cap B$ the free boundary of $u$. A free boundary point $x_0\in \partial\{u\in M\}\cap B$  is said to be {\it regular} if 
$$
\nu_{u(x_0)}\cdot A_{u(x_0)}(Du(x_0),Du(x_0)) < 0,
$$
and 
$$
\limsup_{r\to 0} \frac{\min\diam(\{ u\in\partial M\}\cap B_r(x_0))}{r} > 0.
$$
The point  $x_0\in\partial\{u\in M\}\cap B$ is said to be  {\it singular} if it is {\it not} regular. 
\end{defn}


\section{Optimal regularity for constraint maps}\label{sec:reg}

In this section, we study the optimal regularity of constraint maps. Recall  the notation from   
\eqref{eq:setting}--\eqref{eq:Dut-Dun}.

\begin{lem}\label{lem:ut}
Let $u\in W^{1,2}(\Omega; \overline M)$ be a local constraint map such that $u\in C(B;\cN(\partial  M))$ for some ball $B\subset\Omega$. Then $u\in C^{1,1}_{loc}(B)$, $D^3 V\in BMO_{loc}(B)$, and 
\begin{equation}\label{eq:ut-sys}
\Delta V = -2 (Dw \cdot D)\nu_V+ H_u((D u)^\tau,(D u)^\tau)\quad\text{in }B,
\end{equation}
in the strong sense. 
\end{lem}

\begin{proof}
Since $u\in C(B;\cN(\partial  M))$, Theorem~\ref{thm:prelim} implies that $u\in W_{loc}^{2,p}(B)$ for any $p\in(1,\infty)$, and so it solves \eqref{eq:main}
in the strong sense.
Now, with $J_u : = \nabla \Pi\circ u$ and $H_u := \nabla^2 \Pi \circ u$,
it follows that
$$
\Delta V=J_u \Delta u + H_u (Du, Du).
$$
Since $ \Delta u = A_u(Du,Du)\chi_{\{u\in \partial M\}} = F \nu_u \chi_{\{u\in \partial M\}}$ with $F:= -\nu_u \cdot A_u(Du,Du)$, and $ (\nu\cdot \nabla) \Pi \equiv 0$, we have 
$$
J_u \Delta u = F \chi_{\{u\in \partial M\}} J_u \nu_u = 0, 
$$
hence
\begin{equation}\label{eq:ut-sys-re}
\Delta V= H_u (D u,D u)\quad\text{in }B.
\end{equation}
To derive the right-hand side of \eqref{eq:ut-sys}, we use \eqref{eq:Du-decom}, \eqref{eq:Dut-Dun}, and $\nabla^2\Pi(\nu,\nu) = 0$, which lead us to
$$
\begin{aligned}
 H_u (D u,D u)  =  2 (D_\alpha w) H_u ((D_\alpha u)^\tau, \nu_V)+ H_u((D u)^\tau, (D u)^\tau) .
\end{aligned}
$$
Now, denoting by $\xi^\tau$ the tangential component of $\xi$ with respect to the tangent hyperplane to $\partial M$, it follows from $\nu\equiv \nu\circ \Pi$ and $(\nu\cdot\nabla)\Pi\equiv 0$ that $\nabla^2 \Pi (\xi^\tau,\nu) = -(\nabla\nu) \xi$. Therefore
$$
H_u ((Du)^\tau,\nu_V) = -((\nabla\nu)\circ V) (D u)= - D(\nu\circ V), 
$$
which proves \eqref{eq:ut-sys}. 

Now, if $u\in C(B)$, then by the partial regularity Theorem \ref{thm:prelim} we know that $u\in W_{loc}^{2,p}(B)$ for any finite $p\in(1,\infty)$, hence $H_u(Du,Du) \in W^{1,p}_{loc}(B)$.
By \eqref{eq:ut-sys}, this implies that $V \in W^{3,p}_{loc}(B)$ for any $p<\infty$, and so in particular $V \in C^{1,1}_{loc}(B)$. Recalling that also $w \in C^{1,1}_{loc}(B)$ (see Theorem~\ref{thm:prelim}), we deduce that 
$u = V + w \nu_V \in C^{1,1}_{loc}(B)$.
In particular, this implies that $H_u(Du,Du) \in C_{loc}^{0,1}(B)$, so  $D^3V\in BMO_{loc}(B)$ follows from \eqref{eq:ut-sys-re} and elliptic regularity theory.
\end{proof}

We can also obtain the basic structure of the set of singular free boundary points from the literature on the scalar obstacle problems. 

\begin{lem}\label{lem:FB-sing}
Let $u$ and $B$  be as in Lemma \ref{lem:ut}, and assume that $\nu_V\cdot A_V((Du)^\tau,(Du)^\tau) < 0$ in $B$. Define  $\Sigma$ to be  the set of all singular free boundary points in $B$.\footnote{See Definition \ref{def:FB}.} Then $\Sigma$ is locally contained in a $(n-1)$-dimensional $C^1$-manifold. 
\end{lem}

\begin{proof}
By Lemma \ref{lem:ut} and the Sobolev embedding, $V\in C_{loc}^{2,\sigma}(B)$ for any $\sigma\in(0,1)$, thus $(D u)^\tau \in C_{loc}^{1,\sigma}(B)$ (recall \eqref{eq:Dut-Dun}). This implies that $g :=- \nu_V\cdot A_V((Du)^\tau,(Du)^\tau) \in C_{loc}^{1,\sigma}(B)$. Now, by assumption, $g> 0$ in $B$. So, in view of \eqref{eq:ur-pde}, our conclusion follows from \cite[Theorem 8]{Caff-obs}. 
\end{proof}

We are now ready to prove the optimal regularity of  solutions to our problem.

\begin{proof}[Proof of Theorem \ref{thm:C11}]
By Theorem \ref{thm:prelim}, it suffices to prove that $u\in C^{1,1}(B)$ for any neighborhood $B\Subset\Omega$ where $u\in C(B)$. Then the result follows from Lemma \ref{lem:ut}.
\end{proof}


\section{Higher regularity of free boundaries}\label{sec:higher}

Here we establish the higher-order regularity of free boundaries around their regular points. We shall continue to use the notation introduced in \eqref{eq:setting}--\eqref{eq:Dut-Dun}. Note that, by Definition \ref{def:FB} and Theorem \ref{thm:prelim}, the free boundary is a $C^1$-graph locally around a regular point. 

\begin{proof}[Proof of Theorem \ref{thm:higher}]

Write $E = \{u \in M\}$  (non-coincidence set) and let $x_0 \in B\cap \partial E$ be a regular free-boundary point. By Theorem \ref{thm:prelim}, $B\cap \partial E$ is a $C^{1,\alpha}$-graph for some ball $B\Subset\Omega$ centered at $x_0$, so by elliptic regularity $ w\in C^2(B\cap \overline E)$.

Fix $\sigma\in(0,1)$. As an induction hypothesis, let us assume that $ V\in C_{loc}^{k+2,\sigma}(B\cap\overline E)\cap C_{loc}^{k+2,\sigma}(B \setminus E)$, for some integer $k\geq 0$; note that the hypothesis is clearly true for $k=0$ due to Lemma \ref{lem:ut}. Arguing as in the proof of Lemma \ref{lem:FB-sing}, $(Du)^\tau \in C_{loc}^{k+1,\sigma}(B\cap\overline E)\cap C_{loc}^{k+1,\sigma}(B \setminus E)$, whence $\nu_V \cdot A_V((Du)^\tau , (Du)^\tau)$ belongs to the same class. However, by \eqref{eq:ur-pde}, $ w \in C^2(B\cap\overline E)$ solves 
\begin{equation}\label{eq:ur-pde-re2}
\begin{dcases}
\Delta  w = - \nu_V \cdot A_V((Du)^\tau , (Du)^\tau)&\text{in }B\cap E,\\
 w = |D w | = 0 &\text{on }B\cap\partial E.
\end{dcases}
\end{equation}
Thus, the partial hodograph-Legendre transformation \cite[Theorem 6.17]{PSU} implies that $\partial E\cap B$ is locally a $C^{k+2,\sigma}$-graph, and $ w \in C_{loc}^{k+2,\sigma}(B\cap\overline E)$. Then by \eqref{eq:decom}, we obtain $u\in C_{loc}^{k+2,\sigma}(B\cap\overline E)\cap C_{loc}^{k+2,\sigma}(B \setminus E)$; here we also used  $u =  V$ in $B \setminus E$.

We stress that the partial hodograph-Legendre transformation does not apply to $V$, as the regularity of $V$ and its derivatives along $B\cap\partial E$ are  {\it a priori} unknown.
 In view of \eqref{eq:ut-sys-re}, $V$ verifies an elliptic, diagonal system whose right-hand side $H_u(Du,Du)$ has a jump across $B\cap\partial E$; note that $H_u(Du,Du)$ involves the first term $-2(Dw\cdot D)\nu_V$ in the right-hand side of \eqref{eq:ut-sys} whose derivatives have the jump. Thus, we resort to the regularity theory for transmission problems. The latter observation, along with Theorem \ref{thm:C11}, 
implies $H_u(Du,Du) \in C_{loc}^{k+1,\sigma}(B\cap\overline E)\cap C_{loc}^{k+1,\sigma}(B \setminus E)\cap C^{0,1}(B)$. Now employing \cite[Theorem 1.1]{zhuge}, we deduce that 
$D_e  V \in C_{loc}^{k+2,\sigma}(B\cap\overline E)\cap C_{loc}^{k+2,\sigma}(B \setminus E)$ for every $e\in\partial B_1$. Since it holds for every unit direction $e$, we derive $ V \in C_{loc}^{k+3,\sigma}(B\cap\overline E)\cap C_{loc}^{k+3,\sigma}(B\setminus E)$. 

Consequently, we verify that the induction hypothesis is now met for $k+1$ in place of $k$. Therefore, one may iterate the argument as much as one desires and obtain that $\partial E \cap B$ is a $C^\infty$-graph, 
and $u\in C^\infty(B\cap\overline E) \cap C^\infty(B\setminus E)$. 
\end{proof}


\section{Regularity of the projected image}\label{sec:proj}

In this section, we shall study the regularity of the projected image of a constraint map. In the first subsection, we shall provide some general aspects of its behavior that are true for any dimension (for the ambient space). In the second subsection, we shall prove its optimal regularity in two dimensions.


\subsection{General analysis in arbitrary dimensions}\label{sec:proj-nd}


The study of the optimal regularity of the projected image boils down to the study of the system  
\begin{equation}\label{eq:v-w-pde} 
\begin{dcases}
\Delta v = a \cdot D w\\
\Delta w = g \chi_{\{w > 0\}}
\end{dcases}
\quad\text{in }B_1,
\end{equation} 
where $a$ is vector-valued, $w, g\geq 0$, and $0\in\partial\{w>0\}$. We shall derive this connection at the end of this subsection, see Remark \ref{rem:proj}.

The break of the continuity of $D^3 v$ occurs only at the free boundary points of $w$, due to the possible jump discontinuity of $D^2 w$. As our analysis is of local character, and translation invariant, we shall also assume 
\begin{equation}\label{eq:w-0}
|v|\leq 1,\quad w\geq 0\quad\text{in }B_1, \quad 0\in\partial \{ w> 0\}.
\end{equation}
As \eqref{eq:v-w-pde} comes from our problem \eqref{eq:ur-pde} and \eqref{eq:ut-sys}, it is natural to assume 
\begin{equation}\label{eq:a-g-C1a}
\| a \|_{C^{1,\sigma}(B_1)} \leq 1,\quad 0\leq 2g\leq 3\quad\text{in }B_1,\quad [ g ]_{C^{0,1}(B_1)} \leq \frac{1}{2},
\end{equation}
for some $\sigma\in(0,1)$. In addition, we suppose that
\begin{equation}\label{eq:a-g}
|a| \leq \lambda \sqrt g\quad\text{on }\{w = 0\}, \qquad (\hbox{for some } \lambda>1) ,
\end{equation}
 which will turn out to be a relation between the coefficients involved in \eqref{eq:ur-pde} and \eqref{eq:ut-sys}. 

\begin{rem}\label{rem:v-w}
The following facts will be used frequently in our analysis. By \eqref{eq:a-g-C1a}, the classical $C^{1,1}$-estimate \cite[Lemmas 3, and 4]{Caff-obs} yields 
\begin{equation}\label{eq:w-C11}
\| w \|_{C^{1,1}(B_{3/4})} \leq c,
\end{equation}
for some dimensional constant $c$. Thus, by \eqref{eq:v-w-pde}, $\| \Delta v\|_{C^{0,1}(B_{3/4})} \leq c$, and thus $v\in W^{3,p}(B_{1/2})$, for any $p\in(1,\infty)$, by the $L^p$-theory \cite[Theorem 9.19]{GT}. As a result, for each unit direction $e\in\partial B_1$, $v_e := D_e v \in W^{2,p}(B_{1/2})$ is a strong solution to 
\begin{equation}\label{eq:ve-pde}
\Delta v_e = a\cdot Dw_e + a_e \cdot Dw\quad\text{in }B_{1/2},
\end{equation}
where $a_e := D_e a$ and $w_e := D_e w$, and by the Sobolev embedding,
\begin{equation}\label{eq:v-W3p}
\| v \|_{C^{2,1-\frac{n}{p}}(B_{1/2})} \leq c_p,\quad\forall p \in(n, \infty),
\end{equation}
where $c_p$ depends only on $n$ and $p$. 
\end{rem}

To carry out uniform estimates for scalar obstacle problems, we need $g$ to be uniformly positive. This is {\it not} necessarily the case here. Nevertheless, the issue can be resolved by means of the additional relation \eqref{eq:a-g} between the coefficients $a$ and $g$; see Lemma \ref{lem:rescale} and Remark \ref{rem:rescale} below.

\begin{lem}\label{lem:rescale}
Let $(v,w)$ solve \eqref{eq:v-w-pde} in $B_1$, and assume \eqref{eq:w-0}, \eqref{eq:a-g-C1a}, and \eqref{eq:a-g}. Then there exists a cubic harmonic polynomial $Q$, such that 
\begin{equation}\label{eq:v-Q-1}
\sup_{B_r} | v - Q | \leq c\lambda r^4,\quad\forall r\in\left( g(0),1\right),
\end{equation}
where $c>1$ depends only on $n$.
\end{lem}

\begin{proof}
Write 
\begin{equation}\label{eq:e}
\e := g(0) \geq 0.
\end{equation}
The argument below holds for $0 \leq \e<\frac{1}{2}$; if $\e \geq\frac{1}{2}$, then \eqref{eq:v-Q-1} becomes trivial by choosing $Q\equiv 0$ (recall that $|v|\leq 1$ in $B_1$). In what follows, $c>1$ is a generic dimensional constant. Recalling that $[g ]_{C^{0,1}(B_1)} \leq \frac{1}{2}$, we have
\begin{equation}\label{eq:Lapw}
\sup_{B_r} |\Delta w| \leq \sup_{B_r} |g| \leq \e + \frac{r}{2} \leq  \frac{3r}{2},\quad\forall r \in (\e,1).
\end{equation}
Since $w\geq 0$ in $B_1$ and $w(0) = |Dw(0)| = 0$, \eqref{eq:Lapw} implies
\begin{equation}\label{eq:Dw-C11}
\sup_{B_r} (w + r |Dw|) \leq c r^3,\quad\forall r\in\left(\e,\frac{1}{2}\right);
\end{equation}
see \cite[Lemmas 3, and 4]{Caff-obs}. However, by \eqref{eq:a-g}, \eqref{eq:w-0}, and \eqref{eq:e},
$$
|a(0)| \leq \lambda\sqrt\e. 
$$
which along with $[a]_{C^{0,1}(B_1)}\leq 1$, \eqref{eq:Dw-C11}, and $\lambda>1$, leads us to 
$$
\sup_{B_r} |\Delta v| \leq (\lambda \sqrt\e + r) \sup_{B_r} |Dw| \leq c\lambda r^2\sqrt r,\quad\forall r\in \left(\e,\frac{1}{2}\right). 
$$
As we assume $|v|\leq 1$ in $B_1$, our conclusion \eqref{eq:v-Q-1} follows from the standard iteration technique \cite[Theorem 3]{Caff}. 
We remark that the approximating polynomial $Q$ can be chosen to be harmonic, as $v$ can be approximated by a harmonic function up to scale $\e$.
\end{proof}

\begin{rem}\label{rem:rescale}
Owing to the above lemma (note that $Q$ in \eqref{eq:v-Q-1} is a harmonic polynomial), if $\e := g(0) > 0$, one may consider rescalings
$$
v_\e (y) := \frac{ (v-Q)(\e y)}{c\lambda \e^4},\quad w_\e (y):= \frac{w(\e y)}{\e^3},
$$
which now solve
$$
\begin{dcases}
\Delta v_\e = a_\e \cdot Dw_\e \\
\Delta w_\e = g_\e \chi_{\{w_\e >0\}}
\end{dcases}
\quad\text{in }B_1, 
$$
where $a_\e := a(\e y)/c\lambda $ and $g_\e := g(\e y)/{\e}$. 
Due to \eqref{eq:a-g-C1a}, $g_\e (0) = 1$ and $1\leq 2g_\e \leq 3$ in $B_1$, while by \eqref{eq:v-Q-1}, $(v_\e,w_\e)$ and $(a_\e,g_\e)$ continue to satisfy \eqref{eq:w-0} and \eqref{eq:a-g-C1a}.
\end{rem}

By Remark \ref{rem:rescale}, we may now assume
\begin{equation}\label{eq:a-g-re}
g(0) = 1,
\end{equation}
in addition to \eqref{eq:a-g-C1a}. The additional relation \eqref{eq:a-g} will no longer be  needed from this point. In the two  subsequent  lemmas, we shall consider two special cases 
that ensure cubic growth of $v$ up to a quadratic polynomial.

Our first scenario is when the minimal diameter of $\{w =0\}$ has geometric  decay.
\begin{lem}\label{lem:below}
Let $(v,w)$ solve \eqref{eq:v-w-pde}, and assume \eqref{eq:a-g-C1a}, \eqref{eq:w-0}, and \eqref{eq:a-g-re}. Fix $\gamma\in(0,1)$ and $\lambda>1$, and suppose that there exists some $r_0\in [0,1)$ such that
\begin{equation}\label{eq:below}
\min\diam(\{w = 0\}\cap B_r) \leq \lambda r^{1+\gamma},\quad \forall r\in(r_0,1).
\end{equation}
Then there is a cubic harmonic polynomial $Q$ for which
\begin{equation}\label{eq:v-Q}
\| v - Q \|_{C^1(B_r)}^* \leq c \lambda^{\frac{1}{2n}} r^{3+\bar\gamma},\quad \forall r\in (r_0,1),
\end{equation}
where both $c> 1$ and $\bar\gamma\in(0,\gamma)$ depend only on $n$ and $\gamma$.
\end{lem}

\begin{proof}
Throughout the proof, we shall denote by $c$ a large positive constant depending only on $n$ and $\gamma$. 
We may assume without loss of generality that $r_0 \leq \frac{1}{4}$, since otherwise the conclusion \eqref{eq:v-Q} becomes trivial by choosing $Q\equiv 0$ (recall that $|v|\leq 1$ in $B_1$, see \eqref{eq:w-0}). Let $e\in\partial B_1$ be arbitrary, $v_e$, $w_e$ and $a_e$ be as in Remark \ref{rem:v-w}, and recall \eqref{eq:w-C11}, \eqref{eq:ve-pde}, and \eqref{eq:v-W3p}. 

Recall also that we assume $r_0 \leq\frac{1}{4}$. For each $r\in(r_0,\frac{1}{2})$, let $\psi_r$ solve 
\begin{equation}\label{eq:hr-pde}
\begin{dcases}
\Delta \psi_r = g &\text{in }B_r, \\
\psi_r = w &\text{on }\partial B_r. 
\end{dcases}
\end{equation}
By \eqref{eq:a-g-C1a}, we know that $\| \Delta \psi_r \|_{C^{0,1}(B_r)} \leq \frac{1}{2}$. Due to \eqref{eq:w-0} and \eqref{eq:w-C11}, 
\begin{equation}\label{eq:w-C11-r}
\sup_{B_r} (w + r |Dw|) \leq cr^2.
\end{equation}
Thus, by the Schauder theory \cite[Theorem 6.17]{GT}, 
there is a quadratic polynomial $q_r$ such that 
\begin{equation}\label{eq:hr-C2a}
\| D^2 (\psi_r - q_r) \|_{L^\infty(B_s)} \leq c \left(\frac{s}{r}\right)^\gamma,\quad\forall s\in (0,r).  
\end{equation}
In comparison of \eqref{eq:ur-pde} with \eqref{eq:hr-pde}, $w - \psi_r \in W_0^{2,n}(B_r)$ verifies 
$$
\Delta (w - \psi_r) =-g\chi_{\{w  =0 \}}\quad\text{in }B_r.
$$ 
The global $W^{2,2n}$-estimates \cite[Corollary 9.9]{GT}, along with \eqref{eq:below}, yield that 
\begin{equation}\label{eq:ur-hr-W2n}
\| D^2 (w - \psi_r) \|_{L^p(B_r)}\leq  c |\{w = 0\}\cap B_r|^{\frac{1}{2n}}  
 \leq cr (\lambda r^\gamma)^{\frac{1}{2n}}. 
\end{equation} 
Combining \eqref{eq:hr-C2a} with \eqref{eq:ur-hr-W2n}, we obtain 
\begin{equation}\label{eq:ur-qr-W2n}
\| D^2 (w - q_r) \|_{L^{2n}(B_s)} \leq cr (\lambda r^\gamma)^{\frac{1}{2n}}+ c s \left(\frac{s}{r}\right)^\gamma,\quad \forall s\in(0,r). 
\end{equation}
Hence, choosing $\e := \frac{\gamma}{2n(1+\gamma)}$, $\bar\gamma:= \frac{\e\gamma}{1+\e}$ and taking $r = s^{\frac{1}{1+\e}}$ in \eqref{eq:ur-qr-W2n}, 
\begin{equation}\label{eq:ur-qhs-W2n}
\| D^2 (w - q_r) \|_{L^{2n}(B_{r^{1+\e}})} \leq  c\lambda^{\frac{1}{2n}} r^{(1+\e)(1+\bar\gamma)}.
\end{equation}
Since \eqref{eq:ur-qhs-W2n} holds for any $r\in(r_0,\frac{1}{2})$, and $c$ being independent of $r$, we now obtain, via a simple iteration, a fixed quadratic polynomial $q$ (which may depend on $r_0$) such that 
 \begin{equation}\label{eq:ur-q-W2n}
\| D^2 (w - q) \|_{L^{2n}(B_r)} \leq c \lambda^{\frac{1}{2n}} r^{1+\bar\gamma},\quad\forall r\in \left(r_0^{1+\e},\frac{1}{2}\right);
\end{equation}
we remark that, although $q$ may change as $r_0$ changes, the constant $c$ in \eqref{eq:ur-q-W2n} stays uniformly bounded. Note that $\bar\gamma$ depends only on $n$ and $\gamma$. Moreover, by \eqref{eq:ur-q-W2n}, \eqref{eq:w-C11-r}, and $\lambda>1$,
\begin{equation}\label{eq:q-C11}
|D^2 q| \leq c\lambda^{\frac{1}{2n}}. 
\end{equation}
Now let $v'$ be the solution to  
\begin{equation}\label{eq:vphe-pde}
\begin{dcases}
\Delta v' = a\cdot Dq_e &\text{in } B_{1/2},\\
v' = v_e &\text{on }\partial B_{1/2},
\end{dcases}
\end{equation}
where $q_e := D_e q$. By \eqref{eq:a-g-C1a}, \eqref{eq:q-C11}, and \eqref{eq:v-W3p}, $\| v' \|_{C^{2,\frac{1}{2}}(B_{1/4})} \leq c\lambda^{1/(2n)}$. This immediately implies 
\begin{equation}\label{eq:hve-C3}
\| v' - q' \|_{C^1(B_r)}^* \leq c\lambda^{\frac{1}{2n}} r^{2+\frac{1}{2}},\quad\forall r\in\left(0,\frac{1}{2}\right),
\end{equation}
where $q'$ is the second-order Taylor polynomial of $v'$ at the origin, and $\| \cdot \|_{C^1}^*$ is the adimensional $C^1$-norm, see \eqref{eq:adim}. 

By \eqref{eq:a-g-C1a}, \eqref{eq:ve-pde}, \eqref{eq:w-C11-r}, \eqref{eq:ur-q-W2n}, and \eqref{eq:vphe-pde}, we observe that $v'' := v_e - v' \in W_0^{2,2n}(B_{1/2})$ and 
\begin{equation}\label{eq:D2he-Ln}
\| \Delta v'' \|_{L^{2n}(B_r)} \leq c\lambda^{\frac{1}{2n}} r^{1+\bar\gamma} + cr \sup_{B_r} |Dw| \leq c \lambda^{\frac{1}{2n}} r^{1+\bar\gamma}.
\end{equation}
Since \eqref{eq:D2he-Ln} holds for any $r\in(r_0,\frac{1}{2})$, by standard iteration techniques \cite[Theorem 3]{Caff} we can find a harmonic quadratic polynomial $q''$ such that 
\begin{equation}\label{eq:he-C2a}
\| v'' - q'' \|_{C^1(B_r)}^* \leq c\lambda^{\frac{1}{2n}} r^{2+\bar\gamma},\quad \forall r\in\left(r_0,\frac{1}{2}\right).
\end{equation}
Set $Q^e := q' + q''$. Note that \eqref{eq:v-W3p} implies $D_j v_{e_i} = v_{e_ie_j} = D_i v_{e_j}$ in $B_{1/2}$. Hence, by \eqref{eq:hve-C3} and \eqref{eq:he-C2a}, we have 
$$
\max_{i\neq j}\sup_{B_r} | D_j Q^{e_i} - D_i Q^{e_j} | \leq c \lambda^{\frac{1}{2n}} r^{1+\bar\gamma},\quad\forall r\geq r_0.
$$
Thus, we can find a cubic polynomial $Q$ such that $|D_e Q - Q^e| \leq c\lambda^{1/(2n)} r^{1+\bar\gamma}$ in $B_r$ for all $r > r_0$, for every $e\in \partial B_1$. Then our conclusion \eqref{eq:v-Q} follows from \eqref{eq:hve-C3} and \eqref{eq:he-C2a}. 
\end{proof} 

Another favourable  scenario is when $w - p$ grows subquadratically at infinity after rescaling, for some quadratic polynomial $p$. In this case, we use an estimate on the generalized Newtonian potential to obtain a uniform cubic growth of $v$ up to a quadratic polynomial.

\begin{lem}\label{lem:C21-gen}
Let $(v,w)$ solve \eqref{eq:v-w-pde} in $B_1$, and suppose that \eqref{eq:w-0} and \eqref{eq:a-g-C1a} hold. Let $\delta\in(0,1)$, $\alpha \in[0,1]$, $\lambda>1$ be given,  and assume that there is a quadratic polynomial $p$ for which $|D^2 p| \leq \lambda$ and
\begin{equation}\label{eq:w-p-gen}
\sup_{B_r} | D(w-p) | \leq \lambda
\begin{dcases}
 \delta^\alpha r^{1-\alpha},&\text{if } \alpha\in(0,1], \\
\frac{ r}{|\log\delta|}, &\text{if } \alpha = 0,
\end{dcases}
\quad\forall r\in(\delta,1).
\end{equation}
Then 
\begin{equation}\label{eq:C21-gen}
\sup_{B_r} |v - q| \leq c \lambda r^3,\quad\forall r\in(\delta,1),
\end{equation}
where  $q$  is the second-order Taylor polynomial of $v$ at the origin,  
and  $c$ depends only on $n$, $\alpha$ and $\sigma$. 
\end{lem}

\begin{proof} 
For each unit direction $e\in\partial B_1$, set 
\begin{equation}\label{eq:fi}
f_e := D_e(w-p).
\end{equation}
Then, since $|D^2 p|\leq \lambda$ and $\lambda>1$ (see \eqref{eq:a-g}), recalling \eqref{eq:w-C11} we have $|Df_e| \leq c\lambda$ in $B_{3/4}$. Also, thanks to \eqref{eq:w-p-gen}, $\frac{1}{c\lambda}f_e$ verifies \eqref{eq:f-grow} with $\omega(r) = r^{1-\alpha}$ if $\alpha > 0$ and $\omega(r) \equiv 1/|\log \delta|$ if $\alpha = 0$. Therefore, we can apply Lemma \ref{lem:quad} (with $B_1$ there replaced by $B_{3/4}$) to each component of the generalized (vector-valued) Newtonian potential 
\begin{equation}\label{eq:Vij}
\Phi_e (x) := \int_{B_{3/4}} G(x,y) Df_e (y) \,dy, 
\end{equation}
with $G$ as in \eqref{eq:Hess-fund}. This yields  
\begin{equation}\label{eq:Phi-l} 
 \sup_{B_r} |\Phi_e|  \leq c \lambda r^2,\quad\forall r \in (\delta,1),
\end{equation}
On the other hand, by \eqref{eq:v-W3p}, one can find the second-order Taylor polynomial $q$ for $v$ at the origin, such that 
\begin{equation}\label{eq:v-C11}
\| q \|_{L^\infty(B_1)} \leq c.
\end{equation}
Recall from Remark \ref{rem:v-w} that $v_e := D_e v$ solves \eqref{eq:ve-pde} in the strong sense. Write $q_e := D_e q$. Then by \eqref{eq:fi}, \eqref{eq:Vij}, and \eqref{eq:ve-pde}, the function 
$$
h_e := v_e - q_e - a\cdot \Phi_e,
$$
satisfies 
$$
\Delta h_e = a \cdot D p_e + a_e \cdot Dw\quad \text{in } B_{1/2}, 
$$
in the strong sense. Also by \eqref{eq:a-g-C1a}, \eqref{eq:v-W3p}, \eqref{eq:Phi-l}, and \eqref{eq:v-C11}, we have 
\begin{equation}\label{eq:hi-Linf}
\| h_e \|_{L^\infty(B_{1/2})} \leq c\lambda. 
\end{equation}
Furthermore, by \eqref{eq:a-g-C1a}, \eqref{eq:w-C11}, and $|D^2 p|\leq \lambda$, it holds $\| \Delta h_e \|_{C^{0,1}(B_{1/2})} \leq c\lambda$. So
 the interior $C^{1,1}$-estimate applies to the above PDE, which along with \eqref{eq:hi-Linf} implies that 
\begin{equation}\label{eq:D2hi-Linf}
\| D^2 h_e  \|_{L^\infty(B_{1/4})} \leq c\lambda. 
\end{equation}
By the choice of $q$ and $\Phi_e$, we have $h_e(0) = |Dh_e (0)| = 0$, so  \eqref{eq:D2hi-Linf} leads to
\begin{equation}\label{eq:hd-quad}
\sup_{B_r} |h_e| \leq c \lambda r^2,\quad\forall r\in \left(0,\frac{1}{4}\right).
\end{equation}
By \eqref{eq:Phi-l} and \eqref{eq:hd-quad}, we arrive at 
$$
\sup_{B_r} |v_e - q_e|  \leq  c\lambda r^2,\quad\forall r\in \left(\delta,\frac{1}{4}\right).
$$
As $c$ is independent of the unit direction $e$, and $D^k v(0) = D^k q(0)$ for $k \in \{0,1,2\}$, we obtain \eqref{eq:C21-gen} for $r\in(\delta,\frac{1}{4})$. Thanks to \eqref{eq:v-C11}, one can extend this estimate to all $r\in(\delta,1)$ by enlarging the constant $c$ if necessary. 
\end{proof}

Our last scenario is when the free boundary of $w$ is regular. 

\begin{lem}\label{lem:v-C21-reg}
Let $(v,w)$ be a solution to \eqref{eq:v-w-pde}, and assume \eqref{eq:w-0}, \eqref{eq:a-g-C1a} and \eqref{eq:a-g-re}. Suppose further that $B_1\cap \partial\{ w >0\}$ is a $C^1$-graph. Then $v\in C^{2,1}(B_{1/4})$ and 
\begin{equation}\label{eq:v-C21-reg}
\| v \|_{C^{2,1}(B_{1/4})} \leq c, 
\end{equation}
where $c$ depends only on $n$, $\sigma$ and $\min\diam(\{ w = 0\}\cap B_1)$. 
\end{lem}

\begin{proof}
By \eqref{eq:a-g-C1a}, \eqref{eq:a-g-re}, and \cite[Theorem 7]{Caff-obs}, $B_{3/4}\cap \partial\{w >0\}$ is in fact a $C^{1,\alpha}$-graph, for some $\alpha\in(0,1)$ depending only on $n$. Now for each $e\in\partial B_1$, $w_e := D_e w$ verifies 
$$
\begin{dcases}
\Delta w_e = g_e &\text{in }B_{3/4}\cap\{w>0\},\\
w_e = 0 &\text{on }B_{3/4}\cap\partial\{w>0\},
\end{dcases}
$$
so $w_e \in C^{1,\alpha}(B_{1/2}\cap \overline{\{w>0\}})$, whose norm depends on $n$, $\alpha$, and the $C^{1,\alpha}$-character of $B_{3/4}\cap \partial\{w>0\}$; the last parameter is fully determined by $\min\diam(\{w = 0\}\cap B_1)$.
Now, due to \eqref{eq:w-0}, \eqref{eq:a-g-C1a}, and \eqref{eq:ve-pde}, one can apply \cite[Theorem 1.2]{Xio} to $v_e := D_e v$, for each $e\in\partial B_1$, and obtain 
$$
\| v_e \|_{C^{1,1}(B_{1/4})} \leq c,
$$ 
where $c$ depends only on $n$, $\sigma$ and the $C^{1,\alpha}$-character of $B_{1/2}\cap \partial\{w >0 \}$. Since the estimate holds uniformly for all $e\in\partial B_1$, our conclusion follows. 
\end{proof}

We remark as follows the connection between the projected image of a constraint map and the system \eqref{eq:v-w-pde} for $(v,w)$. 

\begin{rem}\label{rem:proj}
Let $\Omega$ be a bounded domain in $\R^n$, $n\geq 2$, and let $M$ be an open set in $\R^m$, $m\geq 2$, with smooth boundary $\partial M$. Let $u\in W^{1,2}(\Omega;\overline M)$ be a constraint map, such that $u\in C(B; \cN(\partial  M))$ for a ball $B\subset\Omega$, with $\cN(\partial  M)$ the tubular neighborhood. Here $\rho$, $\nu$, $\Pi$ are as in \eqref{eq:setting}, and they are well-defined and smooth. Write $w := \rho \circ u$ and $V:= \Pi\circ u$. 

By Lemma \ref{lem:ut} and \eqref{eq:ut-sys}, $V \in C(\overline B)\cap W^{3,p}(B)$ is a solution to 
$$
\Delta V^i = a^i \cdot D w + h^i, \quad 1\leq i\leq m,
$$
where $a^i := - 2 D(\nu^i\circ u)$, and $h^i:= H^i(u)((Du)^\tau, (Du)^\tau)$. However, by Lemma \ref{lem:ut} and \eqref{eq:Du-decom}, $(Du)^\tau \in C^{1,\sigma}(B)$, so $h^i \in C^{1,\sigma}(B)$. Hence, subtracting from $V^i$ the solution $\vp^i\in C^{3,\sigma}(B)\cap C(\overline B)$ to
$$
\begin{dcases}
\Delta \vp^i = h^i &\text{in }B,\\
\vp^i = V^i&\text{on }\partial B,
\end{dcases}
$$
we observe that $(v,w) := (V^i - \vp^i,w)$ solves \eqref{eq:v-w-pde} in $B$, with $a = a^i$, and $g = -\nu_u\cdot A_u((Du)^\tau,(Du)^\tau)$. Again by $(Du)^\tau\in C^{1,\sigma}(B)$, $g \in C^{1,\sigma}(B)$. The basic assumptions in \eqref{eq:w-0} and \eqref{eq:a-g-C1a} can be verified by suitable rescaling.

To see \eqref{eq:a-g}, one needs to make use of convexity of $\partial M$ on the contact set. Indeed, by the definition of the second fundamental form, for each $\xi\in\R^m$,
$$
\nu\cdot A(\xi^\tau,\xi^\tau) = - \Hess \rho(\xi,\xi)\quad\text{on }\partial M,
$$
where $\xi^\tau$ is the orthogonal projection of $\xi$ onto the tangential hyperplane to $\partial M$. Thanks to Remark~\ref{rem:prelim}, $\Hess\rho$ is a nonnegative definite, symmetric, matrix-valued map on the contact set, $\{u\in \partial M\}$ so one can write it as $\Hess \rho \equiv\beta^2$, for some symmetric matrix-valued $\beta\geq 0$. Since $\nu = \grad \rho$, we also have $\nabla \nu =  \beta^2$. Then $a^i = -2 D(\nu^i \circ V) = -2( (\beta_{ij}^2)\circ V) DV^j = -2 \beta_{ik}\circ V (DV)_k^\beta$, where we wrote $(DV)_k^\beta := (\beta_{kj}\circ V) DV^j$. With $\beta$ at hand, one may also rewrite $g$ as $g \equiv -\nu_V\cdot A_V((DV)^\tau,(DV)^\tau) =  \beta^2(V) (DV,DV) = |(DV)^\beta|^2$. Thus, 
$$
|a^i| \leq 2\left(\sup_\partial M |\beta|\right) \sqrt g\quad\text{on }\{u\in\partial M\},
$$
which implies \eqref{eq:a-g} with $\lambda =2 \sup_\partial M |\beta|$, since $\{w = 0\} = \{u\in\partial M\}$.
\end{rem}

Let us close this subsection with the proofs of Theorems \ref{thm:proj-reg} and \ref{thm:proj-sing}. 

\begin{proof}[Proof of Theorem \ref{thm:proj-reg}]
We already proved that $D^3V\in BMO_{loc}(B)$ in Lemma \ref{lem:ut}. Now suppose that $\partial\{u\in M\}\cap B$ consists of regular points only. Then following the notation in Remark \ref{rem:proj}, the pair $(v^i,w) := (V^i - \vp^i,w)$, $1\leq i\leq m$, verifies the assumption of Lemma \ref{lem:v-C21-reg}. Hence, $D^3(V^i - \vp^i) \in L_{loc}^\infty(B)$, which together with $\vp^i \in C_{loc}^{3,\sigma}(B)$ proves $D^3 V^i\in L_{loc}^\infty(B)$. This finishes the proof. 
\end{proof}


\begin{proof}[Proof of Theorem \ref{thm:proj-sing}]
Let $x_0\in\partial\{u \in M\}\cap B$ be a singular point,
and let $V^i$, $\vp^i$, $w$, and $g$ be as in Remark \ref{rem:proj}. 
In what follows  $c,C$ are positive constants that may depend on $x_0$ but not on $r$. 

If $g(x_0)=0$, then the result follows from Lemma \ref{lem:rescale}. So we can assume that $g(x_0)>0$.
Due to \cite{Caff-obs77}, $w\in C_{loc}^{2,\bar\sigma}(\{x_0\})$ implies 
$$
\sup_{B_r(x_0)} |w - p_{x_0}| \leq Cr^{2+\bar\sigma},\quad\forall r \in (0,r_{x_0}),
$$
for some $r_{x_0} > 0$ and some convex quadratic polynomial $p_{x_0}$ with $\Delta p_{x_0} = g(x_0) > 0$. In particular, the minimal diameter of $\{ p_{x_0} = 0\}$ is equal to zero. Now, if $y_0 \in \{w = 0 \}\cap B_r(x_0)$, then by the quadratic behavior of $p_{x_0}$ it follows that $\dist(y_0,\{p_{x_0} = 0\})^2 \leq c p_{x_0}(y_0) \leq Cr^{2+\bar\sigma}$. This implies that $\{w = 0\}\cap B_r(x_0)\subset \cN_{Cr^{1+\frac{\bar\sigma}{2}}}(\{p_{x_0} = 0\})$.  Thus, 
$$
\min\diam (\{w = 0\}\cap B_r(x_0)) \leq C r^{1+\frac{\bar\sigma}{2}},\quad\forall r\in(0,r_{x_0}).
$$
We can then apply Lemma \ref{lem:below} to $V^i - \vp^i$, which proves that $V^i - \vp^i \in C^{3,\sigma}(\{x_0\})$, with $\sigma \in (0,\bar\sigma)$ depending only on $n$ and $\bar\sigma$. As $\vp^i \in C_{loc}^{3,\sigma}(B)$, we conclude that $V^i \in C^{3,\sigma}(\{x_0\})$. 
\end{proof}



\subsection{Optimal regularity in two dimensions}\label{sec:proj-2d}

The aim of this subsection is to prove Theorem \ref{thm:C21}. In view of Remark \ref{rem:v-w} and \ref{rem:rescale}, we shall continue to assume that $(v,w)$ solves \eqref{eq:v-w-pde} and \eqref{eq:w-0}, \eqref{eq:a-g-C1a}, \eqref{eq:a-g-re} hold, and most importantly, the ambient space is of dimension $n = 2$. In connection with \eqref{eq:v-w-pde}, \eqref{eq:w-0}, and \eqref{eq:a-g-re}, we shall call $U$ a global solution if $U$ is a strong solution to
\begin{equation}\label{eq:U-pde}
\begin{dcases}
\Delta U = \chi_{\{U>0\}}\quad\text{in }\R^2, \\
U\geq 0,\, 0\in\partial\{U>0\}. 
\end{dcases}
\end{equation}
We shall keep this setting throughout this section, unless stated otherwise. From now on, the constant $c$ will vary upon each occurrence, but it will remain to be dependent at most on $n$ and $\sigma$, unless stated otherwise. Also $\| \cdot \|_{C^1}^*$ is the adimensional $C^1$-norm, i.e., 
$$
\| f \|_{C^1(B_r)}^* := \sup_{B_r} |f| + r \sup_{B_r} |Df|. 
$$ 

Using  standard perturbation arguments, Theorem \ref{thm:C21} becomes a simple corollary to the following proposition.

\begin{prop}\label{prop:C21}
Let $v$, $w$, $a$, $g$, and $\sigma$ be as in \eqref{eq:v-w-pde}, \eqref{eq:w-0}, \eqref{eq:a-g-C1a}, and \eqref{eq:a-g-re}. Then there exists $\e_\sigma \in (0,1)$, depending only on $\sigma$, such that the following holds for any $\e\in(0,\e_\sigma)$: if for each $r\in(0,1]$ there exists a global solution $U^r$ to \eqref{eq:U-pde} such that 
\begin{equation}\label{eq:w-Ur}
\| w - U^r \|_{C^1(B_r)}^* \leq \e r^{2+\sigma},\quad [g]_{C^{0,\sigma}(B_1)} \leq \e, 
\end{equation}
then there exists a quadratic polynomial $q$ such that 
\begin{equation}\label{eq:v-q}
\sup_{B_r} |v - q| \leq cr^3,\quad\forall r\in[0,1],
\end{equation}
where $c$ depends at most on $\sigma$. 
\end{prop}

To prove this proposition, the starting point is Lemma \ref{lem:below}, which asserts the following: if the minimal diameter of the contact set $\{ w = 0\}$ decays with order $(1+\gamma)$ up to some scale, say $r_0$, then the solution $v$ to \eqref{eq:v-w-pde} admits $(3+\bar\gamma)$-order approximation by a cubic polynomial up to scale $r_0$, where $\bar\gamma$ depends only on $n$ and $\gamma$. We may take $\gamma = \sigma/4$ with $\sigma$ as in \eqref{eq:w-Ur}. 

If $r_0 = 0$, then we are done. This corresponds to the case that the origin is a singular point of $\partial\{w > 0\}$. 

If $r_0 > 0$, this means that the minimal diameter of $\{w = 0\}$ no longer decays with order $(1+\gamma)$ below the scale $r_0$, but instead it starts to open up. This corresponds to the case when the origin is a regular point of $\partial \{w > 0\}$. Note that $r_0$ is by no means universal, as it depends on the $C^1$-character of $\partial \{w > 0\}$ at the origin.

The following lemma shows how to treat the latter scenario. This is precisely the place where we use the approximation \eqref{eq:w-Ur} by global solutions.

\begin{lem}\label{lem:above}
There exists $\e_\sigma\in(0,1)$, depending only on $\sigma$, such that the following holds whenever $\e < \e_\sigma$: if 
\begin{equation}\label{eq:w-U-g-re}
\| w - U \|_{C^1(B_1)} \leq \e,\quad [g ]_{C^{0,\sigma}(B_1)} \leq \e, 
\end{equation}
and 
\begin{equation}\label{eq:md-w}
\min\diam(\{ w = 0\}\cap  B_1) \geq 2\left(\frac{\e}{\e_\sigma}\right)^{\frac{1}{4}},
\end{equation}
then there exists a quadratic polynomial $q$ such that
$$
\sup_{B_r} |v-q| \leq cr^3,\quad\forall r\in(0,1),
$$
where $c$ depends only on $\sigma$. 
\end{lem}

\begin{proof}
Throughout this proof, $c$ will be a positive constant, possibly depending on $\sigma$, and it may vary upon each occurrence. Thanks to \eqref{eq:w-U-g-re} and \eqref{eq:md-w}, one can compare the free boundary of $w$ to that of the global solution $U$, see \cite[Lemma 12]{Caff-obs}. More precisely, one can find a sufficiently small constant $\e_\sigma\in(0,1)$ such that 
\begin{equation}\label{eq:w-U-FB}
\partial\{w >0\}\cap  B_1 \subset \cN_{c\sqrt\e} (\partial\{ U > 0\}),
\end{equation}
where $\cN_\e(S)$ is the $\e$-neighborhood of set $S$, and 
\begin{equation}\label{eq:md-U-re}
\min\diam(\{ U =0\} \cap  B_1) >\left(\frac{\e}{\e_\sigma}\right)^{\frac{1}{4}}.
\end{equation}
With \eqref{eq:md-U-re} at hand, one may choose 
\begin{equation}\label{eq:d-re}
\delta\in\left[\sqrt{\frac{\e}{\e_\sigma}},1\right]
\end{equation}
such that Lemma \ref{lem:U-p} holds. 

\setcounter{step}{0}
\begin{step}
Estimate in $[\delta,1]$.
\end{step}

By the first assertion \eqref{eq:U-p} of Lemma \ref{lem:U-p}, there exists a quadratic polynomial $p$ with $|\Delta p|\leq 1$ for which 
\begin{equation}\label{eq:U-p-re}
\sup_{B_r} |D(U- p)| \leq c\sqrt{\delta r},\quad\forall r\in[\delta,1].
\end{equation}
Combining \eqref{eq:U-p-re} with \eqref{eq:w-U-g-re}, \eqref{eq:d-re}, and the triangle inequality, yields
\begin{equation}\label{eq:w-p}
\sup_{B_r} |D(w - p)| \leq \e + c\sqrt{\delta r} \leq c\sqrt{\delta r},\quad\forall r\in[\delta,1].
\end{equation}
By \eqref{eq:w-C11} and \eqref{eq:w-p}, we may invoke Lemma \ref{lem:C21-gen} to derive that 
\begin{equation}\label{eq:v-q-re}
\sup_{B_r} |v- q| \leq \hat cr^3,\quad\forall r\in\left[\delta,\frac{2}{3}\right],
\end{equation}
where $q$ is the second-order Taylor polynomial of $v$ at the origin. As $|v|\leq 1$ in $B_1$ and $|a|\leq 1$, \eqref{eq:v-w-pde} and \eqref{eq:w-C11} imply $\|\Delta v\|_{C^{0,1}(B_{3/4})} \leq c$. Therefore, we know that $|q|\leq \hat c$ in $B_1$. One can also enlarge $\hat c$ in \eqref{eq:v-q-re} to extend this estimate to all $r\in (\delta,1)$. 

\begin{step}
Estimate in $(0,\delta)$.
\end{step}

By the second assertion of Lemma \ref{lem:U-p} with $\delta$ as in \eqref{eq:d-re}, one may consider the following two cases only. 

\setcounter{case}{0}

\begin{case}
$\min\diam(\{U = 0\}\cap  B_\delta) \geq \delta/c_0$.
\end{case}

Let us rescale our solutions as 
$$
v_\delta (y) := \frac{(v-q)(\delta y)}{\hat c\delta^3},\quad w_\delta(y) := \frac{w(\delta y)}{\delta^2},
$$
where $\hat c>1$ is as in \eqref{eq:v-q-re}.
As $q$ is the second-order Taylor polynomial of $v$ at $0\in\partial\{w>0\}$, $\Delta q = \Delta v(0) = a(0)\cdot Dw(0) = 0$. Hence, $(v_\delta,w_\delta)$ solves
$$
\begin{dcases}
\Delta v_\delta = a_\delta \cdot Dw_\delta\\
\Delta w_\delta = g_\delta \chi_{\{w_\delta > 0\}}
\end{dcases}
\quad\text{in }B_1,
$$
with $a_\delta (y) := a(\delta y)/c$ and $g_\delta(y) := g(\delta y)$. With these definitions, $(v_\delta,w_\delta)$ and $(a_\delta,g_\delta)$ continue to satisfy \eqref{eq:w-0}, \eqref{eq:a-g-C1a}, and \eqref{eq:a-g-re}. Moreover, by \eqref{eq:w-U-FB} and \eqref{eq:d-re}, we have 
$$
\partial\{w_\delta > 0\}\cap  B_\delta \subset \cN_{c\sqrt{\e_\sigma}} (\partial\{ U_\delta > 0\}),
$$
where $U_\delta(y) := U(\delta y)/\delta^2$. Choosing $\e_\sigma$ smaller if necessary, we obtain from $\min\diam(\{ U_\delta =0\}\cap B_\delta)\geq \delta/c_0$ that 
\begin{equation}\label{eq:md-w-d}
\min\diam(\{w_\delta =0\}\cap  B_1) \geq \frac{1}{2c_0}. 
\end{equation}
Also, by \eqref{eq:w-U-g-re} and \eqref{eq:d-re}, one has 
\begin{equation}\label{eq:g-Ca-d}
[g_\delta ]_{C^{0,\sigma}(B_\delta)} \leq \sqrt{\e_\sigma}. 
\end{equation}
Therefore, taking $\e_\sigma$ even smaller, one may deduce from the uniform regularity \cite[Theorem 7]{Caff-obs} for obstacle problems that $\partial\{w_\delta > 0\}\cap  B_{1/2}$ is a $C^{1,\alpha}$-curve, for some $\alpha\in(0,\sigma)$ depending at most on $n$ and $\sigma$. Now, by Lemma \ref{lem:v-C21-reg}, $v_\delta \in C^{2,1}(B_{1/4})$ and 
\begin{equation}\label{eq:vd-C21}
\| v_\delta\|_{C^{2,1}(B_{1/4})} \leq c,
\end{equation}
where $c$ depends at most on $n$, $\sigma$ and $\min\diam(\{w_\delta =0\}\cap B_1)$; however, thanks to \eqref{eq:md-w-d}, the constant $c$ in \eqref{eq:vd-C21} only depends on $n$ and $\sigma$. 

Rescaling back and recalling again that $q$ is the second-order Taylor polynomial of $v$ at the origin, it follows from \eqref{eq:vd-C21} that 
\begin{equation}\label{eq:v-q-re2}
\sup_{B_r} |v-q| \leq cr^3,\quad\forall r\in \left(0,\frac{\delta}{4}\right). 
\end{equation}
This estimate can easily be extended to $[\frac{\delta}{4},\delta)$ by enlarging $c$, as \eqref{eq:v-q-re} again implies $|v-q|\leq c\delta^3$ in $B_\delta$. Hence, the proof is finished for the case $\min\diam(\{ U =0\}\cap  B_\delta) \geq \delta/c_0$. 

\begin{case}
 $B_\delta\setminus \{U =0\}$ is a disjoint union of two open connected components satisfying \eqref{eq:dist-d} and \eqref{eq:md-re}. 
 \end{case}

Recall \eqref{eq:w-U-FB} again, and let $v_\delta$ and $w_\delta$ be as above. By \eqref{eq:dist-d} (with $\mu^2 = \frac{\e}{\e_\sigma}$ there), we can also choose $\e_\sigma$ small enough such that $B_1\setminus \{w_\delta =0\}$ becomes a disjoint union of two open and connected components $\omega_i$, $i\in\{1,2\}$. Also, by \eqref{eq:w-U-FB} and \eqref{eq:md-re}, each $\omega_i$ satisfies 
\begin{equation}\label{eq:md-re-2}
\min\diam(B_1\setminus\omega_i) \geq \frac{1 - \sqrt{\e_\sigma}}{c} \geq \frac{1}{2c_0}.
\end{equation}
Since $\omega_1\cap\omega_2 = \emptyset$, we can split $w_\delta$ into the sum of $w_\delta^i$ in $B_1$, where $w_\delta^i\geq 0$ in $B_1$ and $\omega_i = B_1\cap \{ w_\delta^i > 0\}$. Then we set $v_\delta^i$ in such a way that $(v_\delta^i,w_\delta^i)$ is a solution to \eqref{eq:v-w-pde}, and $v_\delta^1 = v_\delta$ and $v_\delta^2 = 0$ on $\partial B_1$. In this way, we can repeat the argument in Case 1 to each pair $(v_\delta^i,w_\delta^i)$. The conclusion then follow by combining together the resulting estimates for each pair. 
\end{proof}

We are now ready to establish the proof for the main proposition. 

\begin{proof}[Proof of Proposition \ref{prop:C21}] 
Let $\e_\sigma$ be as in Lemma \ref{lem:above} and define
\begin{equation}\label{eq:r0}
r_0 := \inf\biggl\{ r \in (0,1] : \min\diam (\{ w =0\}\cap B_r ) \leq 2r\left(\frac{\e r^\sigma}{\e_\sigma}\right)^{\frac{1}{4}} \biggr\}. 
\end{equation}
Note that $r_0 \in [0,1]$. In what follows, we will use $c$ to denote a generic constant depending at most on $\sigma$. 

\setcounter{step}{0}

\begin{step}
Estimate in $(r_0,1)$.  
\end{step}

If $r_0 \in [ \frac{3}{4},1]$, one may skip this step. Otherwise, we can apply Lemma \ref{lem:below} with $\lambda = 2\e_\sigma^{-1/4}$ and $\gamma = \sigma/4$. This yields a cubic harmonic polynomial $Q$ satisfying \eqref{eq:v-Q}. Writing by $q$ the quadratic part (i.e., polynomial up to second-order) of $Q$, one deduces from \eqref{eq:v-Q} and the triangle inequality that 
\begin{equation}\label{eq:v-q-below}
\sup_{B_r} |v - q| \leq \hat c r^3,\quad\forall r\in\left(r_0,\frac{3}{4}\right].
\end{equation}
By enlarging $\hat c$, we can extend this estimate to $r\in(\frac{3}{4},1)$. We remark that $q$ is also a harmonic polynomial, as so is $Q$. Let us also remark that if $r_0 = 0$, then our proof is finished.

\begin{step}
Estimate in $(0,r_0)$.
\end{step}
Let $\hat c>1$ be as in \eqref{eq:v-q-below}, and 
rescale our solutions as  
\begin{equation}\label{eq:w0}
v_0 (y) := \frac{(v-q)(r_0 y)}{\hat cr_0^3},\quad w_0 (y) := \frac{w(r_0 y)}{r_0^2},\quad U_0 (y):= \frac{U^{r_0}(r_0y)}{r_0^2},
\end{equation}
where $U^{r_0}$ is the global solution to \eqref{eq:U-pde} for which \eqref{eq:w-Ur} holds at $r=r_0$. As $q$ is a harmonic polynomial and $\hat c>1$, we can argue  as in the proof of Lemma \ref{lem:above},  to conclude  that
 $(v_0,w_0)$ solves \eqref{eq:v-w-pde} with coefficients $(a_0,g_0)$, and the structure conditions \eqref{eq:w-0}, \eqref{eq:a-g-C1a} and \eqref{eq:a-g-re} continue to hold. Moreover,
 \eqref{eq:w-U-g-re} follows directly from 
\eqref{eq:w-Ur}, and  \eqref{eq:md-w} holds by the definition of $r_0$ in \eqref{eq:md-w}. Since $\e r_0^\sigma< \e_\sigma$, with $\e_\sigma$ as in Lemma \ref{lem:below}, we may now apply the lemma to $v_0$ and obtain 
\begin{equation}\label{eq:v0-q0}
\sup_{B_s} |v_0 - q_0| \leq c s^3,\quad\forall s\in (0,1],
\end{equation}
where $q_0$ is the second-order Taylor polynomial of $v_0$ at the origin. Rescaling back from \eqref{eq:v0-q0} and combining it with \eqref{eq:v-q-below}, yields the desired estimate \eqref{eq:v-q}. 
\end{proof}

We are ready to prove the optimal regularity of the image map in two dimensions.

\begin{proof}[Proof of Theorem \ref{thm:C21}]
Let $\Omega\subset \R^2$ be a bounded domain and  $u\in W^{1,2}(\Omega;\overline M)$ be a constraint map. By the partial regularity \cite[Theorem 1.1]{DF} of constraint maps in two-dimension, and Theorem \ref{thm:C11}, $u\in C_{loc}^{1,1}(\Omega)$. Let $\{\rho,\nu,\Pi\}$ be as in \eqref{eq:setting}, and set $w := \rho\circ u$ and $V := \Pi\circ u$. 

Let $x_0\in\partial\{u\in M\}\cap\Omega$ be given, and let $B\Subset\Omega$ be the ball centered at $x_0$ such that $\diam B = d_{x_0} := \dist(x_0,\partial\Omega)$. Fix $i\in\{1,2,\cdots,m\}$ and let $a^i$, $h^i$, $g$ and $\vp^i$ be as in Remark \ref{rem:proj}, and set $v^i := V^i - \vp^i$. Let $c>1$ be a large constant determined solely by $m$, the $C^{3,\sigma}$-character of $\partial M$, $\diam\Omega$, and $\diam u(\Omega)$. Then the quantities $\|a^i\|_{C^{1,\sigma}(B)}^*$, $\|h^i\|_{C^{1,\sigma}(B)}^*$ and $\|g \|_{C^{0,1}(B)}^*$ can all be bounded by $c$; the definition of the adimensional norm $\|\cdot\|^*$ can be found in the beginning of Section \ref{sec:prelim}. 

As $\Delta \vp^i = h^i \in C^{1,\sigma}(B)$ and $\vp^i = V^i$ on $\partial B$, Schauder theory yields
\begin{equation}\label{eq:vpi-C21}
\sup_{rB} | \vp^i - \hat q_{x_0}^i | \leq c\left(\frac{r}{d_{x_0}}\right)^3,\quad\forall r\in(0,1).
\end{equation}
In addition, \eqref{eq:v-w-pde}, \eqref{eq:w-0}, \eqref{eq:a-g-C1a}, and \eqref{eq:a-g} are satisfied by suitably scaled versions of $(v^i,w)$ and $(a^i,g)$; the scaling should also map $B$ to the unit disk $B_1$. At this point, one can follow the procedure of Lemma \ref{lem:rescale}, Remark \ref{rem:rescale}, and Proposition \ref{prop:C21}, to obtain a quadratic polynomial $q_{x_0}^i$ such that 
\begin{equation}\label{eq:vi-C21}
\sup_{rB} |v^i - q_{x_0}^i | \leq cr^3,\quad\forall r\in(0,1).
\end{equation}
Here the estimate is independent of the size of $g(x_0)$, due to Lemma \ref{lem:rescale}. If $g(x_0) = 0$, this lemma yields the full estimate. Otherwise, i.e., if $g(x_0) > 0$, after additional rescaling shown in Remark \ref{rem:rescale}, we may apply Proposition \ref{prop:C21} by the uniform geometric approximation property of $w$, see Definition \ref{def:gap}. We also remark that \eqref{eq:vi-C21} is independent of $d_{x_0}^{-3}$, as $|v^i| \leq d_{x_0}^3$ in $B$; this is because $v^i = 0$ on $\partial B$ and $|\Delta v^i| \leq |a^i||Dw| \leq c d_{x_0}$ in $B$, where the last inequality is deduced from \eqref{eq:w-C11}. 

Combining \eqref{eq:vpi-C21} and \eqref{eq:vi-C21}, since $V^i = v^i + \vp^i$ we obtain the desired estimate, namely interior cubic growth of $V^i$ up to a quadratic polynomial at every free boundary point $x_0\in\partial\{u\in M\}\cap \Omega$. From here, standard techniques imply that $V^i \in C_{loc}^{2,1}(\Omega)$ for each $i\in\{1,2,\cdots,m\}$. We leave out the details to the reader.
\end{proof}


\appendix

 
 \section{Generalization to ``leaky''   maps}\label{app:gen}
 
Although \eqref{eq:main} can be derived as the Euler-Lagrange system for constraint maps, the system itself may admit non-variational solutions. Moreover, non-variational solution maps do not need to favour one side of $\partial M$ over the other; this can be done by considering the signed distance to $\partial M$ (i.e., extending the distance function $\rho$ negatively to the complement of $\overline M$). Below we describe some possible examples of non-variational solution maps to our system.

\begin{rem}\label{rem:ex}
Let $\Omega$ be a domain in the ambient space $\R^n$, say $0\in\partial \Omega$, $B$ be a ball around the origin, and $u : B\setminus \Omega$ be a continuous harmonic map; we remark that the continuity of $u$ already imposes some topological condition on the portion $u(B\setminus \Omega)$ of $\partial M$. Now suppose that $u$ extends to $B$ so that $\Delta u =0$ 
in $B\cap\Omega$ and $u\in C^1(\Omega)$; i.e., not only the values but also the gradient matches on the boundary $B\cap\partial \Omega$. Under such  matching conditions, it is not difficult to verify that the extended version solves our system \eqref{eq:main} in the sense of distribution, see \eqref{eq:weak}. Note that $u|_{B\cap\Omega}$ no longer needs to favour one side of $\partial M$ to another; in fact, $u$ may ``leak'' from one side of $\partial M$ to another. Our main interest here is to answer the following question: ``What can be said for both $u$ and $\Omega$ when such an extension holds?" This problem naturally extends the connection  between the scalar no-sign obstacle problem and the harmonic continuation property of a given domain, see \cite{PSU}. 
\end{rem}

Let $M\subset\R^m$ be a smooth domain with boundary $\partial M$. Then there is a tubular neighborhood $\cN(\partial  M)$ where the decomposition \eqref{eq:setting} is valid; here we extend $\rho$ to $-\dist(\cdot,\partial M)$ in $M^c$. Let $B\subset\R^n$ be a ball. 
 
We shall call $u\in W^{1,2}(B;\cN(\partial  M))$ a solution to \eqref{eq:main} in the distributional sense, if 
\begin{equation}\label{eq:weak}
\begin{aligned}
\int_B D u : D \vp\,dx& = - \int_{B\cap\{ u\in\partial M\}} A_u(Du, Du) \cdot \vp \,dx ,
\end{aligned}
\end{equation}
for any $\vp\in C_0^\infty(B;\R^m)$.

\begin{lem}\label{lem:prelim}
Let $u\in C\cap W^{1,2}(B;\cN(\partial  M))$ be a solution to \eqref{eq:main} in the distributional sense. Then $u\in W_{loc}^{2,p}(B)$ for any $1 \leq p < \infty$.
\end{lem}

\begin{proof}
By the regularity theory for quadratic systems \cite[Theorem 9.7, 9.8]{GM}, $u\in C_{loc}^{1,\sigma}(B)$ for all $\sigma\in(0,1)$. In particular, $|Du| \in L_{loc}^\infty(B)$, which in turn implies that $|\Delta u| \in L_{loc}^\infty(B)$. The conclusion now follows immediately from the classical $L^p$-theory. 
\end{proof}

Let us now use the notations in \eqref{eq:setting} -- \eqref{eq:Dut-Dun}. By Lemma \ref{lem:prelim}, we may now work with strong solutions $u\in W^{2,p}(B)$ to \eqref{eq:main}, $p > n$. 

\begin{lem}\label{lem:ur}
Let $u\in W^{2,p}(B;\cN(\partial  M))$ be a strong solution to \eqref{eq:main}. Then $w \in C_{loc}^{1,1}(B)$, and it solves 
\begin{equation}\label{eq:ur-pde-re}
\Delta  w=  - \nu_V\cdot A_V((Du)^\tau,(Du)^\tau) \chi_{\{  w \neq 0\}}
\end{equation}
in the strong sense. 
\end{lem}

\begin{proof}
By the chain rule and \eqref{eq:setting}, 
$$
\begin{aligned}
\Delta  w &=  \grad\rho(u) \cdot \Delta u + \Hess\rho(u)(Du,Du)\\
&= \nu_u\cdot A_u(Du,Du)\chi_{\{ u \in\partial M\}} +  \Hess\rho(u)(Du,Du) = - \nu_V\cdot A_V((Du)^\tau,(Du)^\tau) \chi_{\{  w \neq 0\}},
\end{aligned}
$$
a.e.\ in $B$, where the last identity follows from $\Hess\rho(\xi,\xi) = - \nu\cdot A(\xi^\tau,\xi^\tau)$, with $\xi^\tau$ being the orthogonal projection to the tangent hyperplane to $\partial M$. 

Since $u\in W^{2,p}(B)$, then $(Du)^\tau \in W^{1,p}(B)$ by \eqref{eq:Dut-Dun}. Therefore, since $p>n$, $ \nu_V\cdot A_V((Du)^\tau,(Du)^\tau)  \in C_{loc}^{0,\sigma}(\Omega)$ with $\sigma = 1 - \frac{n}{p}$ (by the Sobolev embedding). According to \cite[Theorem 1.2]{ALS}, $ w \in C_{loc}^{1,1}(B)$. 
\end{proof} 

\begin{lem}\label{lem:ut-re}
Let $u\in W^{2,p}(B;\cN(\partial  M))$ be a strong solution to \eqref{eq:main}. Then $V\in C_{loc}^{3,p}(B)$ and it solves \eqref{eq:ut-sys} in the classical sense. 
\end{lem}

\begin{proof}
Since $u\in W^{2,p}(B)$, $(\Delta u)^\tau \equiv (I - \nu_V\otimes \nu_V) \Delta u = 0$ a.e.\ in $B$. Noticing that $\nu\equiv \nu\circ \Pi$ and $(\nu\cdot\nabla)\Pi \equiv 0$ in $\cN(\partial  M)$, one derives \eqref{eq:ut-sys-re} from the chain rule. The rest of the proof is then the same as in the proof of Lemma \ref{lem:ut}. 
\end{proof}

As a result of Lemmas \ref{lem:prelim} -- \ref{lem:ut-re}, we obtain the optimal regularity of continuous ``leaky'' maps $u$. 

\begin{thm}\label{thm:C11-leak}
Let $u\in C(B;\cN(\partial  M))$ be a solution to \eqref{eq:main} in the distributional sense. Then $u\in C_{loc}^{1,1}(B)$. 
\end{thm}

Some remarks follow.

\begin{rem}\label{rem:leak}
To study the regularity of the free boundary, one needs to put additional assumption that $\nu_V \cdot A_V((Du)^\tau, (Du)^\tau) \leq 0$ in $B$. While for constraint maps this follows immediately from the constraint, but this may not be true in general for leaky maps. Of course, leaky maps satisfying such a sign-condition are those of interests. In fact, the terminology ``leaky map'' makes sense only in this case, as according to the literature on scalar no-sign obstacle problems, $w := \rho\circ u$ is always nonnegative around regular free boundary points. 
\end{rem}

One may notice that the proof for the higher regularity of the free boundaries of constraint maps, Theorem \ref{thm:higher}, does not need the side condition. Thus, we can immediately extend the result to continuous leaky maps. 

\begin{thm} 
Let $u \in C\cap W^{1,2}(\Omega;\cN(\partial M))$ be a solution to \eqref{eq:main} in the distributional sense, and $x_0 \in \partial \{u\not\in \partial M\}\cap \Omega$ be a regular point. Then there is a ball $B\subset\Omega$ centered at $x_0$ such that $\partial\{ u \not\in \partial M\}\cap B$ is a $C^\infty$-graph, and $u\in C^\infty(B\cap \overline{\{ u \not\in \partial M\}})\cap C^\infty(B\cap \{u \in \partial M\}))$. 
\end{thm}

Analogously, we can also extend some of regularity results on the projected image. 

\begin{thm}
Let $B$ be a ball in $\R^n$, $n\geq 2$, and $u\in C\cap W^{1,2}(B;\cN(\partial M))$ be a solution to \eqref{eq:main} in the distributional sense. Then $D^3(\Pi\circ u)\in BMO_{loc}(B)$, and if every point $x_0\in \partial\{u\in M\}\cap B$ is regular then $D^3(\Pi\circ u) \in L_{loc}^\infty(B)$. 
\end{thm}

To extend our result in dimension two (Theorem \ref{thm:C21}), one needs extra care due to the lack of the one-side-condition, especially for the proofs of Lemma 6.2 and Remark 6.3.
In addition, to run the argument based on the geometric approximating property, one needs to understand how to deal with global solutions to the no-sign obstacle problem with wild behavior \cite{Shap1}. Hence, our current argument applies only near points where 
$\nu_V\cdot A_V((Du)^\tau, (Du)^\tau)<0$  (and so, Theorem \ref{thm:C21} extends to weak solutions of \eqref{eq:main} around such points), but otherwise a much more refined analysis is needed.


 
 \section{Properties of Schwarz function}\label{app:schwarz}

Here we establish a uniform subquadratic growth of $U - p$ for any global solution in the plane. We only prove it in dimension $n=2$, as our proof relies on the explicit formula \cite{S} for the (gradient of) global solutions via the Schwarz functions. Still, we believe such a result to hold also in higher dimension.

By \cite{S}, the contact set of global solutions in the plane can only be one of the following:
\begin{itemize}
\item A set enclosed by either an ellipse or a parabola;
\item A strip;
\item A half-plane;
\item A line.
\end{itemize}
Let us identify $\R^2$ with $\C$, and write $z = x + i y$. Since $U$ is a global solution to \eqref{eq:U-pde} in the plane, one can explicitly express the derivative of $U$, according to \cite{S}, as
\begin{equation}\label{eq:DU}
\frac{\partial U}{\partial z} = \frac{1}{4} (\bar z - S(z))\quad\text{in }\{U > 0\},
\end{equation}
where $S$ is the Schwarz function of the analytic curve $\partial\{U > 0\}$; see the monograph \cite{Davis} for more on the Schwarz functions.

\begin{lem}\label{lem:U-p}
Let $n= 2$, and $U$ be a global solution to \eqref{eq:U-pde} and  $\mu := \frac{1}{2}\min\diam(\{U =0\}\cap B_1) > 0$. There is a constant $\delta \in [\mu^2,1]$, and an absolute constant $c_0$, such that the following hold: \medskip

\begin{enumerate}[(i)]
\item {\bf Outside $B_\delta$:} There is a quadratic polynomial $p$, with $\Delta p =  1$,  such that  
\begin{equation}\label{eq:U-p}
\sup_{B_r} |D(U - p)| \leq c_0 \sqrt{\delta r},\quad\forall r\in[\delta,1].
\end{equation}

\medskip

\item {\bf Inside $B_\delta$:} Either 
\begin{equation}\label{eq:md}
\min\diam(\{U =0\}\cap B_\delta) \geq \frac{\delta}{c_0},
\end{equation}
or $\{U > 0 \}\cap B_\delta = \Omega_1\cup \Omega_2$ for some disjoint pair of open connected components $\Omega_i$, $i\in\{1,2\}$, such that the Hausdorff distance $d_H$ satisfies 
\begin{equation}\label{eq:dist-d}
\dist_H (\partial\Omega_1\cap B_\delta, \partial\Omega_2\cap B_\delta) \geq \frac{\mu\sqrt\delta}{c_0}, 
\end{equation}
and
\begin{equation}\label{eq:md-re}
\min\diam(B_\delta \setminus \Omega_i) \geq \frac{\delta}{c_0}\qquad \text{for each $i\in\{1,2\}$.}
\end{equation}
\end{enumerate}
\end{lem}

\begin{proof} 

If $\{U = 0\}$ is either a line, a strip or half-plane, the assertion is immediate. Hence, we shall consider the case where $\{U = 0\}$ is enclosed by either a parabola or an ellipse.

One may write, after rotating the coordinate system, 
\begin{equation}\label{eq:U}
\partial\{ U > 0\} = \{ (x,y) : a^2 x^2 + y^2  = \alpha x + \beta y\},
\end{equation}
where $0\leq a\leq 1$, $\alpha > 0$ and $\beta\in\R$. From now on, we shall use squares $Q_r = (-\frac{r}{2},\frac{r}{2})^2$, instead of disks $B_r$, in order to simplify the exposition.  Moreover, by symmetry, we can assume  that $\beta \geq 0$.
Also, we can assume that $\min\diam(\{U = 0\}\cap Q_2) < 1$ and $Q_2\setminus \{ U =0\}$ is a single connected component, otherwise assertion (ii) of the lemma holds with $\delta = 1$. 

Note that if $a = 0$ then $\partial\{U>0\}$ is a parabola, and $a>0$ corresponds to the case $\partial\{U > 0\}$ being an ellipse. However, even if $a > 0$, depending on its value with respect to $\alpha$, $\partial\{U > 0\}$ may look like a parabola in the unit disk $B_1$. Thus, one has to be careful when studying the latter case.

Since $a\leq 1$ and $\min\diam(\{U = 0\}\cap Q_2) < 1$, we have
\begin{equation}\label{eq:mu}
\begin{aligned}
\mu^2 & := \left[ \frac{1}{2} \min\diam(\{ U =0\}\cap Q_2)\right]^2  = \frac{\beta^2}{4} + \alpha \left(\frac{\alpha}{2a^2} \wedge 1\right) - a^2 \left(\frac{\alpha}{2a^2} \wedge 1\right)^2 < \frac{1}{4}. 
\end{aligned}
\end{equation}
This can be deduced from the fact that $\min\diam(\{U =0\}\cap Q_2)$ is the difference between the roots of $a^2x^2 + y^2 = \alpha x + \beta y$ at $x = \frac{\alpha}{2a^2}\wedge  1$, where $\frac{\alpha}{2a^2}$ is the $x$-coordinate of the center of $\partial\{U>0\}$ when it is an ellipse.

The fact that $Q_2\setminus \{ U =0\}$ is a single connected component, along with $a\leq 1$ and $\alpha > 0$, yields that 
\begin{equation}\label{eq:rho}
\begin{aligned} 
\rho &:= \left| \inf\{ x : (x,y) \in \partial\{U >0\}\} \right| = \frac{\beta^2}{2\alpha}\left[\sqrt{1+ \frac{a^2\beta^2}{\alpha^2}} +1\right]^{-1} < 1.
\end{aligned}
\end{equation} 
This depends on two facts:  first, $\alpha > 0$ implies that $\partial\{U>0\}$ has less portion in $\{(x,y): x\leq 0\}$ than in $\{(x,y): x\geq 0\}$, i.e.,
\begin{equation}\label{eq:portion}
| \inf\{x : (x,y) \in \partial\{ U > 0\}\} | < |\sup\{ x : (x,y) : \partial\{U>0\}\}|;
\end{equation}
second, whether $\partial\{U>0\}$ is an ellipse or a parabola, it has a tip on the leftmost side at $(-\rho,\frac{\beta}{2})$.

We shall divide the proof into two cases, as now the shape of $\partial\{U> 0\}$ becomes more important. This is the key step, and $n = 2$ plays a crucial role. 

\setcounter{case}{0}

\begin{case}
$\alpha > 2a^2$. 
\end{case}

We  claim that both assertions of this lemma holds with  
\begin{equation}\label{eq:d-para}
\delta := \mu^2 \vee \frac{\rho}{4} = \left(\frac{\beta^2}{4} + \alpha - a^2\right)\vee \frac{\rho}{4},
\end{equation}
where $\mu^2$ and $\rho$ are as in \eqref{eq:mu} and \eqref{eq:rho} respectively. Since $0< \mu^2 < 1$ and $0\leq \rho < 1$, we have $0< \delta < \frac{1}{2}$. Moreover, since $\alpha > 2a^2$ implies in \eqref{eq:mu} that $\beta^2\leq 6\alpha$, since $\mu^2 = \frac{\beta^2}{4} + \alpha - a^2$ we have
\begin{equation}\label{eq:md-para}
\begin{aligned}
[\min\diam(\{U =0\}\cap Q_\delta)]^2 &= 4\left(\frac{\beta^2}{4} + \alpha\delta - a^2 \delta^2\right)
\geq 4\mu^2 \delta;
\end{aligned}
\end{equation}
note that the first identity in \eqref{eq:md-para} follows from the fact that, as $\alpha > 2a^2$, the minimal diameter of $\{U =0 \}\cap Q_\delta$ with $\delta < 1$ is equal to the difference between two roots of $a^2 x^2 + y^2 =\alpha x + \beta y$ when $x = \delta$. 

Suppose that $\mu^2\geq \frac{\rho}{4}$ so that $\delta = \mu^2$. Then it follows immediately from \eqref{eq:md-para} that $\min\diam(\{U =0\}\cap Q_\delta) \geq 2\delta$. This verifies the first alternative \eqref{eq:md} in the second assertion of this lemma. 

Next, let us consider the case $\mu^2 \leq \frac{\rho}{4}$, so that $\delta =\frac{ \rho}{4}$. By \eqref{eq:rho} and \eqref{eq:portion}, $Q_\delta=Q_{\rho/4}$ cannot contain any horizontal tip of $\partial\{ U >0\}$. Thus, there are a disjoint pair of open connected components $\Omega_i$, such that $\{ U > 0\}\cap Q_\delta = \Omega_1\cup\Omega_2$. Because of \eqref{eq:portion} and $\alpha > 2a^2$, the Hausdorff distance between $\partial\Omega_1\cap Q_\delta$ and $\partial\Omega_2\cap Q_\delta$ is attained as the difference of two roots of $a^2 x^2 + y^2 = \alpha x + \beta y$ when $x = -\delta = -\frac{\rho}{4}$. As \eqref{eq:rho} implies $\beta^2> 16\alpha \delta$ and since $0 <\delta<\frac{1}{4}$, we can deduce that 
$$
\begin{aligned}
[\dist_H(\partial\Omega_1\cap Q_\delta, \partial\Omega_2\cap Q_\delta)]^2 &= 4\left( \frac{\beta^2}{4} - \alpha\delta- a^2\delta^2\right) 
\geq 4\mu^2\delta,
\end{aligned}
$$
which now proves \eqref{eq:dist-d}. The last assertion \eqref{eq:md-re} can be easily verified as follows. Since we assume $\delta = \frac{\rho}{4} \geq \mu^2 > 0$, we must have $\beta^2 > 0$ by \eqref{eq:rho}. We may assume without loss of generality that $\beta > 0$, as the argument is symmetric for the other case. Then since $y = \frac{\beta}{2}$ divides equally $Q_\delta\cap \{U =0\}$, we may set $\Omega_1$ as the connected component (on the top) for which $Q_\delta\setminus \Omega_1$ contains the lower half-square, $\{(x,y): y \leq 0\}$. Thus, \eqref{eq:md-re} holds with $\Omega_1$. With the choice of $\Omega_1$, the component $\Omega_2$ now lies below $y < \frac{\beta}{2} < \sqrt{12/13}\delta$; the second inequality can be deduced from \eqref{eq:rho} and \eqref{eq:d-para}. Therefore, \eqref{eq:md-re} also holds for $\Omega_2$, which proves that the second assertion of this lemma is satisfied when $\mu^2 \leq\frac{\rho}{4}$ as well.

Now we need to verify  the first part of this lemma. Let us make use of the explicit formula \eqref{eq:DU} for the gradient of $U$. Keeping in mind of the fact that $\alpha > 2a^2$ implies $a$ small, let us substitute $(x,y)$ with $(\frac{z+ \bar z}{2},\frac{z - \bar z}{2i})$  in \eqref{eq:U}, and solve the resulting equation for $\bar z$. It leads us to
\begin{equation}\label{eq:schwarz}
\bar z = \frac{1 + a}{1-a} z - \frac{2az + \zeta}{1-a^2} + \frac{1}{1-a^2}\sqrt{(2az + \zeta)^2 - 4(1-a)\bar\xi z},
\end{equation}
where $\zeta = \alpha + \beta i$ and $\xi = \alpha + a\beta i$. This is the Schwarz function of $\partial\{U > 0\}$ given as in \eqref{eq:U}.

Choose $p$ as the homogeneous, quadratic polynomial satisfying 
\begin{equation}\label{eq:p-para}
\frac{\partial p}{\partial z} = \frac{1}{4} \left( \bar z - \frac{1+a}{1-a} z\right) = - \frac{1}{2} \left(i\Im z + \frac{a}{1-a}z\right). \footnote{Note that $\frac{\partial}{\partial z} = \frac{1}{2}( \frac{\partial }{\partial x} - i \frac{\partial }{\partial y})$.}
\end{equation}
That is, $p(x,y )= (2(1-a))^{-1}( -a x^2 + y^2)$, which clearly verifies $\Delta p = 1$; it is noteworthy that $p$ is not convex. Note that \eqref{eq:mu} and \eqref{eq:rho} together with \eqref{eq:d-para} and $\alpha > 2a^2$ imply $\frac{\beta^2}{4} + \alpha - a^2 \leq \delta < 1$ and $\beta^2 < 6\alpha$. Using $\beta^2 < 6\alpha$ and $2a^2\leq \alpha$, we deduce from \eqref{eq:rho} and \eqref{eq:d-para} again that $\beta^2 \leq c\alpha\delta\leq c\delta^2$. Putting these altogether, 
\begin{equation}\label{eq:z-x}
a < \frac{1}{2}\wedge (c\sqrt\delta), \quad |\zeta| = \sqrt{\alpha^2 + \beta^2} \leq c\delta,\quad |\xi| = \sqrt{\alpha^2 + a^2\beta^2} \leq c\delta. 
\end{equation} 
Utilizing \eqref{eq:z-x} as well as \eqref{eq:DU}, \eqref{eq:schwarz} and \eqref{eq:p-para}, we may proceed as
\begin{equation}\label{eq:U-p-para}
\begin{aligned}
\sup_{B_r\cap \{ U > 0\}} \left| \frac{\partial (U - p)}{\partial z} \right| &\leq c(a r +|\zeta|) + c\sqrt{a^2r^2 + |\zeta|^2 + 4|\xi|r} \\
&\leq c ( \sqrt\delta r + \delta ) + c \sqrt{\delta r^2 + \delta^2 + \delta r}  \leq c \sqrt{\delta r},
\end{aligned}
\end{equation}
for any $r \in (\delta, 1)$. 

To estimate $|D(U - p)| = |Dp|$ in $\{ U = 0\}$, we observe from \eqref{eq:U} that $(\Im z)^2 \leq \alpha|\Re z| + \beta|\Im z| \leq |\zeta| |z|$ for any $z\in\{U =0\}$. Employing \eqref{eq:z-x} as well, we can infer from the second identity in \eqref{eq:p-para} that
\begin{equation}\label{eq:U-p-para-re}
\sup_{B_r\cap \{U =0\}}\left| \frac{\partial p}{\partial z} \right| \leq \sup_{z\in B_r\cap \{U =0\}} \left( \frac{|\Im z|}{2} + \frac{a}{2(1-a)}|z|\right)  \leq c\sqrt{\delta r},
\end{equation} 
for any $r \in (0,1)$. Since $\partial U / \partial z \equiv 0$ in $\{ U =0\}$, the first part of this lemma is now verified by \eqref{eq:U-p-para} and \eqref{eq:U-p-para-re}. 

\begin{case}
$\alpha \leq 2a^2$. 
\end{case}

As for this case, we shall verify the assertions of this lemma with 
\begin{equation}\label{eq:d-ellip}
\delta :=  \mu = \frac{1}{2}\sqrt{\beta^2 + \frac{\alpha^2}{a^2}} .
\end{equation} 
In view of \eqref{eq:mu}, with $\alpha\leq 2a^2$ at hand, we obtain $\min\diam(\{U = 0\}\cap Q_\delta) = \min\diam(\{U = 0\}\cap Q_1) = \delta$. Therefore, the second assertion of this lemma is verified with $\delta$ as in \eqref{eq:d-ellip}.

Thus, it remains for us to prove the first assertion of the lemma. In view of \eqref{eq:U}, we can rewrite $\partial\{U>0\}$ as 
\begin{equation}\label{eq:U-re}
\partial\{U >0\} = z_0 + \delta \left\{ (x,y) : x^2 + \frac{y^2}{a^2} = 1\right\},
\end{equation}
where $z_0 = \frac{\alpha}{2a^2} +i \frac{\beta}{2}$. One may compute the Schwarz function for $x^2 + (y^2/a^2) = 1$ (note  $a\in[\sqrt{\frac{\alpha}{2}},1]\subset(0,1]$) and use dilation to verify that 
\begin{equation}\label{eq:schwarz-re}
\bar z = \bar z_0 + \frac{1-a}{1+a} (z-z_0) - \frac{2a\delta^2}{z-z_0 + \sqrt{z-z_0 - \delta^2(1-a^2)}}
\end{equation}
is the Schwarz function for $\partial\{U>0\}$, now given as in \eqref{eq:U-re}. Set $p$ as the real-valued quadratic polynomial satisfying 
\begin{equation}\label{eq:p-ellip}
\frac{\partial p}{\partial z} = \frac{1}{4} \left( \overline{z-z_0} - \frac{1-a}{1+a} (z-z_0)\right).
\end{equation}
Again $\Delta p = 1$. However, as $\zeta\mapsto \zeta + \sqrt{\zeta^2 - (1-a^2)}$ is a conformal mapping that maps the exterior of the ellipse $E_a = \{(x,y) : x^2 + (y/a)^2 \leq 1\}$ onto the exterior of the unit disk, 
\begin{equation}\label{eq:U-p-ellip}
\sup_{B_1\cap \{U > 0\}} \left|\frac{\partial (U - p)}{\partial z}\right| \leq \frac{2a\delta}{\inf_{\zeta\not\in E_a} | \zeta + \sqrt{\zeta^2 - (1-a^2)}|} \leq c\delta.
\end{equation}
Moreover, in $\{U = 0\} = z_0 + \delta E_a$, one may compute from \eqref{eq:p-ellip} that 
\begin{equation}\label{eq:U-p-ellip-re}
\sup_{B_1\cap \{U =0\}} \left|\frac{\partial p}{\partial z}\right| \leq \delta \sup_{\zeta\in E_a}\left| \bar\zeta - \frac{1-a}{1+a}\zeta\right| \leq c\delta.
\end{equation}
Combining \eqref{eq:U-p-ellip} and \eqref{eq:U-p-ellip-re}, and noting that $\delta \leq \sqrt{\delta r}$ for any $r\in[\delta,1]$, we verify that the first assertion of this lemma holds with $\delta$ as in \eqref{eq:d-ellip}, for the case $\alpha \leq 2a^2$. This concludes the proof for all cases. 
\end{proof}

\section{An estimate for Newtonian potentials}\label{app:newt}

Let $ n\geq 2$, and $\Gamma$ be the fundamental solution of the Laplacian, i.e., 
$$
\Gamma(x) := \begin{dcases}
\frac{1}{2\pi} \log |x|, &\text{if } n= 2, \\
\frac{1}{n(2-n)\alpha_n} |x|^{2-n},&\text{if }n\geq 3,
\end{dcases}
$$
where $\alpha_n$ is the volume of the $n$-dimensional ball. Define
\begin{equation}\label{eq:Hess-fund}
\begin{aligned}
G(x,y) := \Gamma(x-y) - \Gamma(y) + D\Gamma(y)\cdot x,
\end{aligned}
\end{equation}
Given $f\in L^1(\Omega)$, we shall call $V_f := \int_\Omega G(\cdot,y) f(y)\,dy$ the generalized Newtonian potential of $f$, as $\Delta V = f$ in $\Omega$ whenever $V_f \in L^1(\Omega)$. 

\begin{lem}\label{lem:quad}
Let $\omega:(0,\infty)\to(0,\infty)$ be a nondecreasing function with $\omega(1)\leq 1$, $\delta \in (0,1)$ be a given constant, and $f\in C^{0,1}(B_1)$ be such that 
\begin{equation}\label{eq:f-grow}
\sup_{B_r} |f| \leq r\omega\left(\frac{\delta}{r}\right), \quad \forall r\in(\delta,1),
\end{equation}
and for some $1\leq i\leq n$, 
\begin{equation}\label{eq:Df-Linf}
\sup_{B_1} |D_i f| \leq 1. 
\end{equation}
Let $G$ be as in \eqref{eq:Hess-fund} and define, 
$$
\Phi(x) := \int_{B_1} G (x,y) D_i f(y)\,dy.
$$
Then
\begin{equation}\label{eq:quad}
\sup_{B_r} |\Phi| \leq c r\left( \delta + r\omega\left(\frac{\delta}{r}\right) + r\int_{\delta}^{\frac{\delta}{r}}\frac{\omega(\tau)}{\tau}\,d\tau\right),\quad  \forall r \in (\delta,1),
\end{equation} 
where $c>0$ depends only on $n$ and $\omega$. 
\end{lem}

\begin{proof}

By \eqref{eq:Df-Linf}, it is not difficult to prove that \eqref{eq:quad} holds when $r = 1$. Thus, it suffices to prove this estimate for $r\in(\delta,\frac{1}{4})$, provided that $\delta < \frac{1}{8}$. 

Fix any $x\in B_{1/4}\setminus B_\delta$ and write $r := \delta + |x|$, so that $B_\delta(x)\subset B_{2r} \subset B_1$. Integrating by parts, we observe that
\begin{equation}\label{eq:Phi-IJK}
\begin{aligned}
\Phi(x) &= \int_{B_\delta} G(x,y) D_i f(y)\,dy - \int_{B_1\setminus B_\delta}  \frac{\partial  G(x,y)}{\partial y_i} f(y)\,dy  +  \int_{\partial (B_1\setminus B_\delta)} G(x,y) f(y)\nu_i\,d\sigma_y  \\
 &=  \int_{B_\delta} G(x,y) D_i f(y)\,dy - \int_{B_{2r} \setminus B_\delta}  \frac{\partial G(x,y)}{\partial y_i} f(y) \,dy \\
 &\quad - \int_{B_1\setminus B_{2r}}  \frac{\partial G(x,y)}{\partial y_i} f(y) \,dy + \int_{\partial (B_1\setminus B_\delta)} G(x,y) f(y)\nu_i\,d\sigma_y, 
\end{aligned}
\end{equation} 
where $\nu$ is the outward unit normal to the boundary of $B_1 \setminus B_\delta$.

To estimate the first integral on the third line, we observe from \eqref{eq:Df-Linf} that
\begin{equation}\label{eq:Phi-Bd}
\begin{aligned}
\left| \int_{B_\delta} G(x,y) D_i f(y)\,dy\right| &\leq c |x_j |\int_0^1 \,dt \int_{B_\delta} |D_j \Gamma(y-tx) - D_j\Gamma (y) |\,dy \\
&\leq cr \int_0^1 \,dt\int_{B_\delta} \left(\frac{1}{|y-tx|^{n-1}} + \frac{1}{|y|^{n-1}}\right)dy \leq cr \int_{B_\delta} \frac{dy}{|y|^{n-1}} \leq c\delta r,
\end{aligned}
\end{equation}
where to derive the third inequality we used that $\int_{B_\rho(z)} |x-y|^{1-n}\,dy $ is maximized when $z = x$, for any $\rho>0$.

Next, by \eqref{eq:f-grow}, 
\begin{equation}\label{eq:I2}
\begin{aligned}
\int_{B_{2r}\setminus B_\delta}\left|\frac{\partial \Gamma(x- y)}{\partial y_i}\right|  |f(y)|  dy  &\leq cr\omega\left(\frac{\delta}{r}\right) \int_{B_{2r}} \frac{dy}{|x - y|^{n-1}} \leq c r^2\omega\left(\frac{\delta}{r}\right).
\end{aligned}
\end{equation}
Also the symmetry of the Hessian of $\Gamma$ implies 
$$
\int_{B_\rho\setminus B_\delta} D_{ij}\Gamma(y) \,dy = 0,\quad\forall \rho> \delta,
$$
so \eqref{eq:Df-Linf} implies that 
\begin{equation}\label{eq:DI}
\begin{aligned}
\left| x_i \int_{B_{2r}\setminus B_\delta}D_{ij}\Gamma(y)  f(y)dy \right| &\leq cr \int_{B_{2r}\setminus B_\delta } \frac{|f(y) - f(0)|}{|y|^n}\,dy \leq  cr \int_{B_{2r}} \frac{dy}{|y|^{n-1}} \leq cr^2. 
\end{aligned}
\end{equation}
By \eqref{eq:I2} and \eqref{eq:DI}, we have 
\begin{equation}\label{eq:I-l}
\int_{B_{2r}} \left|\frac{\partial G(x,y)}{\partial y_i}\right| |f(y)| dy \leq cr^2\left( 1 + \omega\left(\frac{\delta}{r}\right)\right). 
\end{equation}
This yields the estimate for the second integral on the third line of \eqref{eq:Phi-IJK}.

To estimate the third integral, we observe that $2|y - tsx| \geq 1 + |y|$ whenever $|y| > r$ and $s,t \in [0,1]$, since $r = 1+|x|$. Hence, 
\begin{equation}\label{eq:J-l}
\begin{aligned}
\left| \int_{B_1\setminus B_{2r}} \frac{\partial G(x,y)}{\partial y_i}  f(y) \,dy\right|  
&\leq |x_k x_\ell |\int_0^1 t\,dt \int_0^1  \,ds\int_{B_1\setminus B_{2r}} |f(y) ||D_{ik\ell} \Gamma(y-tsx)| \,dy \\
&\leq cr^2 \int_{B_1\setminus B_{2r}}\omega\left(\frac{\delta}{|y|}\right)\frac{dy}{|y|^n} \leq cr^2 \int_{\delta}^{\frac{\delta}{r}} \frac{\omega(\tau)}{\tau}\,d\tau.
\end{aligned}
\end{equation}
We may analogously estimate the last boundary integral in \eqref{eq:Phi-IJK}. Since again we have $2|y - tsx| \geq 2$ when $|y| = 1 > 2r$, 
\begin{equation}\label{eq:K-l}
\begin{aligned}
\left| \int_{\partial B_1} G(x,y) f(y)\nu_i\,d\sigma_y \right|
&\leq |x_k x_\ell| \int_0^1 t\,dt \int_0^1 \,ds \int_{\partial B_1} |D_{k\ell}\Gamma(y-tsx) f(y)| \,d\sigma_y \\
&\leq c r^2 \int_{\partial B_1} \omega\left(\frac{\delta}{1+|y|}\right)\frac{d\sigma_y}{(1 + |y|)^{n-1}}  \leq c r^2 \omega(\delta). 
\end{aligned}
\end{equation}
Thanks to \eqref{eq:Phi-IJK}, \eqref{eq:I-l}, \eqref{eq:J-l}, and \eqref{eq:K-l},
our proof is complete.
\end{proof}


\end{document}